\documentclass[11pt]{article}
\usepackage[a4paper,margin=2.5cm]{geometry}
\usepackage{amsmath,amssymb,amsfonts,amsthm}
\usepackage[english]{babel}
\usepackage{lmodern}
\usepackage{graphicx}
\usepackage{xcolor}
\usepackage{etoolbox}
\usepackage{algorithm}
\usepackage{algpseudocode}
\usepackage[font=small,labelfont=bf]{caption}
\usepackage[colorlinks=true,allcolors=blue]{hyperref}

\newtheorem{theorem}{Theorem}
\newtheorem{corollary}{Corollary}
\newtheorem{lemma}{Lemma}

\newtheorem{assumption}{Assumption}

\newtheorem{remark}{Remark}

\makeatletter
\newcommand{\articletype}[1]{}
\renewcommand{\title}[1]{\hypersetup{pdftitle={#1}}{\centering\LARGE\bfseries #1\par}\vspace{1.5em}}
\renewcommand{\author}[1]{{\centering\large #1\par}\vspace{1em}}
\newcommand{\affil}[1]{{\centering\small #1\par}}
\newcommand{\email}[1]{\vspace{0.5em}{\centering\small E-mail: #1\par}\vspace{1em}}
\newcommand{\keywords}[1]{\gdef\@arxivkeywords{#1}}
\AtEndEnvironment{abstract}{\ifdefined\@arxivkeywords
  \par\medskip\noindent\textbf{Keywords:} \@arxivkeywords\fi}
\makeatother

\def\xd{x^{\delta}}
\def\yd{y^{\delta}}
\def\E{\mathbb{E}}
\def\P{\mathbb{P}}
\def\F{\mathcal{F}}
\def\Fp{F^{\prime}}
\def\xdag{x^{\dagger}}
\def\nd{n_{\delta}}
\def\N{\mathbb{N}}

\newcommand{\norm}[1]{\left\Vert #1 \right\Vert}
\newcommand{\normm}[1]{\Vert #1 \Vert}
\newcommand{\Ex}[1]{\E\left[ #1 \right]}
\newcommand{\abv}[1]{\left\vert #1 \right\vert}
\newcommand{\brc}[1]{\left\lbrace #1 \right\rbrace}
\newcommand{\Exx}[1]{\E[ #1 ]}

\newcommand{\domain}[0]{\mathrm{U}}
\newcommand{\DefSigma}[0]{\mathcal D _\domain}
\newcommand{\DiscRatio}[0]{R}
\newcommand{\LLdmnd}[1]{\mathbf L^2_{\diamond}\left( #1 \right)}
\newcommand{\Hs}[2]{\mathbf H^{#1}\left(#2\right)}
\newcommand{\Hdmnd}[2]{\mathbf H^{#1}_{\diamond}\left( #2 \right)}
\newcommand{\Lp}[2]{\mathbf L^{#1}\left(#2\right)}
\newcommand{\sigmad}[0]{\sigma^\delta}

\counterwithin{theorem}{section}
\counterwithin{lemma}{section}
\counterwithin{corollary}{section}
\counterwithin{proposition}{section}
\counterwithin{assumption}{section}
\counterwithin{definition}{section}
\counterwithin{remark}{section}
\counterwithin{problem}{section}

\usepackage{todonotes}
\usepackage{comment}
\usepackage{enumitem}
\usepackage{tikz}
\usepackage{pgfplots}
\usepackage{multicol}
\pgfplotsset{compat=1.18}
\usetikzlibrary{external}
\tikzset{external/export=false} 
\definecolor{wue}{RGB}{37,86,149}
\pgfplotsset{
	colormap={wuemap}{rgb255(0cm)=(255,255,255); rgb255(1cm)=(160,195,225); rgb255(2cm)=(37,86,149); rgb255(3cm)=(15,35,70)}
}

\begin{document}

\articletype{Paper} 

\title{A Batched Discrepancy Principle for Stochastic Landweber–Kaczmarz Iteration and its Application to EIT}

\author{Rommel Rabongue Real$^1$, Simon Schwab$^{2,*}$ and Frank Werner$^3$}

\affil{$^1$Department of Mathematics, Physics and Computer Science, University of the Philippines Mindanao, Davao City, Philippines}

\affil{$^2$ Institut für Mathematik, Julius-Maximilians-Universität Würzburg, Würzburg, Germany}

\affil{$^3$ Institut für Mathematik, Julius-Maximilians-Universität Würzburg, Würzburg, Germany}

\affil{$^*$Author to whom any correspondence should be addressed.}

\email{rrreal@up.edu.ph, simon.schwab@uni-wuerzburg.de, frank.werner@uni-wuerzburg.de}


\begin{abstract}
We consider an inverse problem given as ill-posed system of possibly nonlinear operator equations between Hilbert spaces. To solve the problem we use a \emph{stochastic mini-batch Landweber-Kaczmarz method}, drawing uniformly a subset of equations in each step. As an \emph{a posteriori} stopping rule, we introduce a batched version of the discrepancy principle, in which only the residual of the selected batches is monitored. Under the tangential cone condition, we establish the convergence of the iterates for exact data, the convergence of the regularized iterates for noisy data under an \emph{a priori} choice of the stopping index, and the convergence of the method when it is terminated by the proposed batched discrepancy principle. Numerical experiments for the electrical impedance tomography problem illustrate the theoretical results and demonstrate the performance of the method for different batch sizes.
\end{abstract}

\section{Introduction}\label{Section:intro}
Consider inverse problems in terms of a system of $p$ equations
\begin{equation}\label{eq1}
F_{i}(x)= y_{i}, \quad i =1,\ldots,p,
\end{equation}
 where each $F_{i}:\mathcal{D}(F_{i})\subset X  \rightarrow Y_{i}$ is a possibly nonlinear operator between Hilbert spaces $X$ and $Y_{i}$, each with inner products $\langle \cdot,\cdot \rangle$ and norms $\norm{\cdot}$, respectively, and with domains $\mathcal{D}(F_{i})$. Such problems arise in many physical phenomena, such as in imaging, where the number of equations $p$ can be very large.\\
In practical applications, the exact data $y:=(y_{1},\ldots,y_{p})$ are not available. Instead we only have noisy data $y^{\delta}:=(y^{\delta}_{1},\ldots,y^{\delta}_{p})$ satisfying
\begin{equation}
\Vert y_{i}^{\delta} - y_{i} \Vert \leq \delta_{i},
\end{equation}
with noise levels $\delta_i>0$ for $i=1,\ldots, p$.
Since inverse problems are often ill-posed, small perturbations on the data may lead to large perturbations in the corresponding solution. Approximating the solution of \eqref{eq1} using noisy data $\yd$ is done by regularization, an important area of research, see eq. the monographs \cite{Engl2000,Kaltenbacher2008} and the references therein. The various regularization methods, which have been developed to solve single operator equations, can be applied to \eqref{eq1} by reformulating the problem as 
\begin{equation}\label{eq:fullnl1}
F(x) = y,
\end{equation}
where $F:=(F_1,\ldots,F_p) : \mathcal D(F) \subseteq X \to Y$, $Y=Y_1 \times \ldots \times Y_p$ and $\mathcal D(F) = \cap_{i=1}^{p}\mathcal{D}(F_{i
})$.
Assuming that $F$ has a continuous Fr\'{e}chet derivative $F^{\prime}(\cdot)$, one of the most notable iterative regularization methods is the Landweber iteration, with the form
\begin{equation}\label{landw1}
\xd_{n+1} = \xd_{n} - \Fp(\xd_{n})^{\ast}(F(\xd_n) - \yd),
\end{equation}
where $F^{\prime}(\cdot)^{\ast}$ is the adjoint of the Fr\'{e}chet derivative $F^{\prime}$ of $F$. The convergence of iteration \eqref{landw1} was established in \cite{Hanke1995} when the stopping index $\nd$ is chosen by the discrepancy principle 
\begin{equation}\label{landwDP}
    \norm{F(\xd_{\nd}) - \yd } \leq \tau \delta < \norm{F(\xd_{n}) - \yd}, \quad 0 \leq n < \nd,
\end{equation}
with a constant $\tau>1$  and $\delta:= \sqrt{\delta_1^2 + \cdots + \delta_{p}^2}$. Due to its simple implementation and relatively small complexity per iteration compared to variational methods (such as Tikhonov regularization), Landweber iteration has received prominent attention in the inverse problem community, with generalizations to non-smooth operators \cite{Clason2019} and semi-iterative accelerations \cite{Hanke_1991}. 
For convergence and regularization properties of the iteration we refer to the monographs \cite{Engl2000,Kaltenbacher2008}.
However, a single step in \eqref{landw1} requires evaluating the operators $F_i$ and $(F^\prime_i(x))^\ast$ for all $p$ equations simultaneously, which can make the Landweber iteration inefficient in applications with a large number of equations. 
A way to overcome this disadvantage is to apply Kaczmarz-type methods, which require evaluating only one of the $p$ equations in \eqref{eq1} per iteration.
 The classical Landweber-Kaczmarz iteration for solving \eqref{eq1} is given by
\begin{equation}\label{landwkacz1}
    \xd_{n+1} = \xd_{n} - p \Fp_{[n]}(\xd_{n})^{\ast}(F_{[n]}(\xd_n) - \yd_{[n]} ),
\end{equation}
where $[n]:= n\mod{p}\in \lbrace 0,1,\ldots,p-1 \rbrace$. Notice that each of the $p$ equations in \eqref{eq1} is chosen cyclically per iteration in \eqref{landwkacz1}. The regularizing properties of this method have been studied in \cite{Kowar_2002, Haltmeier2007}. In the related literature, the factor $p$ in \eqref{landwkacz1} is neglected, which corresponds to rescaling the operator $F$. However, here the factor $p$ is considered for consistency with the rest of the work. Generalizations of the method like the Steepest-Descent-Kaczmarz method \cite{de2008steepest} also include changing step sizes. Other variants of the method are the projective Landweber-Kaczmarz method \cite{Leitao2016}, the averaged Kaczmarz method \cite{Li2018},  or a Kaczmarz method with convex penalty functional \cite{Jin2013}. However, we point out that the efficiency of the method \eqref{landwkacz1} depends on how the equations in \eqref{eq1} are ordered \cite{Herman1993}. To address this issue, the randomized Kaczmarz method, introduced by Strohmer and Vershynin \cite{Strohmer2009} for consistent linear systems, chooses the equations in \eqref{landwkacz1} randomly instead of cyclically. Closely related is the stochastic gradient descent (SGD) method, which can be seen as generalization of the randomized Kaczmarz method \cite{needell2014stochastic} and is for our problem defined as
\begin{equation}\label{sgd1}
    \xd_{n+1} = \xd_{n} - p \Fp_{i_n}(\xd_{n})^{\ast}(F_{i_n}(\xd_n) - \yd_{i_n} ),
\end{equation}
where the index $i_n$ is drawn uniformly from the set $\lbrace 1,\ldots,p \rbrace$. SGD was pioneered by Robbins and Monro in statistical inference \cite{Robbins1951}. Recent studies exhibit the attractiveness of SGD in large-scale problems that require efficient and practical methods, such as machine learning \cite{Bottou2010, Bottou2018}, and also inspired subsequent work such as \cite{Gitman2019, Hertrich2022} and the references therein. The authors in \cite{Chen2018} used SGD in optical tomography, exhibiting promising numerical results, while the paper \cite{Jia2010} studies SGD in inverse lithography technology. Meanwhile, for inverse problems of the form \eqref{eq1}, recent and generalized convergence analysis of \eqref{sgd1} under noisy data with an a priori stopping rule can be read in \cite{JinLu2019} for the linear and in \cite{Jin2020} for the nonlinear case. The generalization to Banach spaces is analyzed in \cite{JinKereta2023,Jin2026a}. Different generalizations include stochastic mirror descent \cite{Jin2023b,Huang2025} and SGD with convex penalty \cite{Gu2025}. Note that in the SGD literature the iteration is commonly defined with a possibly changing step size $\eta_n$ instead of the factor $p$.

However, while SGD substantially lowers the computational burden per iteration compared with the full Landweber iteration, relying on a single equation per step can lead to high variance in the update direction and consequently to oscillatory, slowly converging behavior. 
For large-scale optimization problems
various variance reduction methods have been proposed to accelerate the stochastic gradient method, inspired by machine learning problems, e.g. the stochastic variance reduced gradient (SVRG) \cite{Johnson2013}, the stochastic dual coordinate ascent (SDCA) \cite{ShalevShwartz2013}, the stochastic averaged gradient (SAG) \cite{schmidt2017minimizing}, the SAGA algorithm \cite{Defazio2014}, the stochastic recursive gradient algorithm (SARAH) \cite{Nguyen2017}, the epoch mixed gradient descent (EMGD) \cite{Zhang2013} and the stochastic primal-dual coordinate method (SPDC) \cite{Zhang2015}. In particular, for linear inverse problems the SVRG was analyzed in \cite{Jin_2022SVRG,Jin2025c,Jin2025b}. However their approaches involve one full Landweber step similar to  \eqref{landw1}, as well as storing past gradients, and mixing them with full gradients. \\
Here we consider a mini-batch Kaczmarz approach. This natural compromise between the full Landweber iteration \eqref{landw1} and the single-equation SGD scheme selects randomly, at every iteration $n$, a subset $I_{n}=\{i_1,\ldots,i_b\}\subset\{1,\ldots,p\}$ of fixed batch size $b:= |I_n|\in\{1,\ldots,p\}$, and updates the iterate via
\begin{subequations}
\label{xn}
\begin{align}
	\label{xn:iteration}
	x^{\delta}_{n+1} &:= x^{\delta}_{n} - \omega_n\frac{p}{b} \sum_{i \in I_{n}} F_{i}'(x^{\delta}_{n})^{*}\left(F_{i}(x^{\delta}_{n}) - y^{\delta}_{i}\right), 
	\end{align}
with an additional step control $\omega_n$, defined as
   \begin{align}
   \omega_n = \begin{cases}
		1 \quad \sum_{i\in I_n}\norm{F_{i}(\xd_{n}) - \yd_{i}}^2\geq\tau^2\sum_{i\in I_n}\delta_{i}^2,\\
		0 \quad \text{else}.\label{xn:step_control}
	\end{cases}
\end{align}
\end{subequations}
The step control monitors the residual of the selected equations and decides, based on the discrepancy principle, whether or not the proposed step is taken.
\begin{remark}
	To reduce superfluous repetitions in the calculation of the residuum, additional numerical adjustments can be added. If $\omega_n=0$ then the elements of $I_n$ could be excluded in the selection of the next batch $I_{n+1}\subseteq \{1,\ldots,p\}\setminus I_n$. Going on with $\omega_{n+1}=0$ then $I_{n+2}\subseteq \{1,\ldots,p\} \setminus ( I_n \cup I_{n+1})$ and so on. When $w_{n+i}=1$ one starts again with all indices. And when $\{1,\ldots,p\}\setminus ( I_n \cup \ldots \cup I_{n+i})= \emptyset$, the classical discrepancy principle \eqref{landwDP} is fulfilled and the algorithm stops. This corresponds to a stochastic version of the termination criterion introduced\linebreak in \cite{Haltmeier2007} for the deterministic Kaczmarz method. 
\end{remark}
If $b=p$ the iteration \eqref{xn}  recovers the full Landweber iteration \eqref{landw1} until the classical discrepancy principle \eqref{landwDP} holds and remains constant then. For $b=1$ the iteration reduces, if $\omega_n=1$, to the single-equation SGD. The batch size $b$ interpolates between these two extremes, trading per-iteration cost against variance of the update. \\
For this mini-batch iteration \eqref{xn} we prove convergence in expectation for noise-free data and the corresponding regularization property for \emph{a priori} stopping. This extends the analysis of \cite{Jin2023b} to nonlinear operators. Additionally we propose and analyze a novel \emph{a posteriori} stopping rule: The \emph{batched discrepancy principle} monitors the batched residual of the last $L_\delta\in \N$ iterations, and is defined as
\begin{equation}
			\label{eq:summed_batch_DP}
		n_\delta=\underset{n\geq L_\delta-1}{\arg\min}\Bigg\{ \sum_{j=0}^{L_\delta-1}\sum_{i\in I_{n-j}}\norm{F_i (\xd_{n-j}) - \yd_i}^2-\tau^2\delta_i^2\leq0\Bigg\}.
\end{equation}
The parameter $L_\delta$ is chosen a priori and corresponds to the length of the aggregation. In this work we prove under the tangential cone condition that for properly chosen $\tau>1$ it holds $\nd<\infty$ (Lemma \ref{lemma:finite_termination}) and the iterate $\xd_{n_\delta}$ converges for any choice $L_\delta \to \infty$ to a solution of problem \eqref{eq1} as $\delta \to 0$ (Theorem \ref{thm:batch_DP_conv}).\\   
The proposed stopping rule goes beyond the work in \cite{Jin2023b,Huang2025}, which analyze a mini-batch approach of the stochastic mirror descent algorithm. In \cite{Jin2023b} the method is applied on linear inverse problems and convergence is proven only for a priori stopping. Building on that, \cite{Huang2025} considers the same iteration on nonlinear operators and proposes an a posteriori stopping rule involving, differently to our approach, the full residuum $\norm{F(\xd_n) - \yd}$. Further works, that analyze approaches using the full residuum, are \cite{Jahn2020,Gu2025} for SGD and \cite{Jin2025c,Jin2025b} for SVRG. Note that especially for a large number of equations $p$ and small batch sizes $b$ the batched version of the discrepancy principle significantly reduces the computational cost per iteration compared to techniques which consider all equations simultaneously. \\
 The batched discrepancy principle can be seen as an instance of \textit{early stopping}, a technique which received some attention in recent works on (statistical) inverse problems, see e.g. \cite{BlanchardHoffmannReiss2018a,BlanchardHoffmannReiss2018b,HuckerReiss2025}. The rationale behind is that for iterative methods, the iteration number is on the one hand directly determining the computational cost (to be kept as small as possible), and on the other hand acts as a decreasing regularization parameter (to be chosen not too small). For methods such as truncated singular value decomposition or conjugate gradients, it has been shown in the above references, that early stopping reduces the computational effort while maintaining the optimality of the underlying method, both theoretically and numerically. \\
 
The remainder of this paper is organized as follows. In Section \ref{section:convergence} we first establish the general setting of the analysis. Then we proof the convergence of the mini-batch Landweber-Kaczmarz iteration \eqref{xn} in three cases. First with noise-free data, continuing with noisy data, proving convergence under a priori stopping and then under the a posteriori stopping rule defined in equation \eqref{eq:summed_batch_DP}. Finally, in Section \ref{Section:numerics} we illustrate the theoretical findings with numerical experiments for the electrical impedance tomography (EIT) problem, a nonlinear ill-posed system with a large number of measurement equations.

\section{Convergence}
\label{section:convergence}
In order to establish convergence of method \eqref{xn}, we first define the space on which the analysis will be carried out. We follow a similar approach as in \cite{Jin2023b}, but correct the definition of the sample space for the batch size $b$.
Let \begin{equation}
    \label{eq:sample_space}
    \Lambda_{b,p}= \left\{ I \subseteq \{1,\ldots,p\} : |I| =b  \right\} 
\end{equation} be the sample space and let $\mathcal{S}_{b,p}$ denote the $\sigma$-algebra consisting of all subsets of $\Lambda_{b,p}$. At each iteration step, a random element of $\Lambda_{b,p}$ is chosen based on the uniform distribution. This means $\mathbb P_{b,p}(I)=1/\binom p  b$ for all $I\in \Lambda_{b,p}$, 
since there are $\binom p b $ different $b$-sized subsets of $\{1,\ldots,p\}$. Therefore, for each $n \geq 1$, we consider $\xd_{n}$ on the sample space

\[
\Lambda_{b,p}^{n}:=\underbrace{\Lambda_{b,p} \times \cdots \times \Lambda_{b,p}}_{n  \text{ copies}},
\]
equipped with the $\sigma$-algebra $\mathcal{S}_{b,p}^{ n}:= \bigotimes_{k=1}^n \mathcal S_{b,p}$ and the product probability measure denoted by $$\mathbb{P}_{n}((I_0,\ldots,I_{n-1})):=\mathbb P_{b,p}(I_0) \cdot \ldots \cdot \mathbb P_{b,p}(I_{n-1}) \text{  for } (I_0,\ldots,I_{n-1})\in \Lambda_{b,p}^n .$$
According to the Kolmogorov extension theorem \cite{Bhattacharya2016}, there exists a unique probability $\mathbb{P}$ defined on the measurable space $(\Omega,\mathcal{F}):=(\Lambda_{b,p}^{\infty},\mathcal{S}_{b,p}^{\infty})$, such that each $\mathbb{P}_{n}$ is consistent with $\mathbb{P}$. Using this, the stochastic mini-batch Landweber-Kaczmarz iteration \eqref{xn} corresponds to a random variable on the probability space 
\begin{equation}
\label{eq:prob_space}    
(\Omega,\F,\P):=(\Lambda_{b,p}^{\infty},\mathcal{S}_{b,p}^{\infty},\P).
\end{equation}
Note that also the stopping index $n_\delta$ of the batched discrepancy principle \eqref{eq:summed_batch_DP} is a random variable on $(\Omega,\mathcal F,\P)$.
We denote $\E$ as expectation on $(\Omega,\F,\P)$ and define $L^{2}(\Omega,X)$ as the space consisting of all random variables $x$ with values in $X$ such that $\E[\Vert x \Vert^{2}]$ is finite; this is a Hilbert space under the norm
\[
\left( \E[\Vert x \Vert^{2}] \right)^{1/2}.
\]
As a proof tool, we will use $\lbrace \F_{n} \rbrace$ to denote the natural filtration, where $\F_{n} = \sigma\left( I_{0},\ldots,I_{n-1} \right)$.
The conditional expectation $\Exx{\cdot \vert\mathcal F_n}$ averages over all future batch selections $l\geq n$ and yields an $\mathcal F_n$-measurable function, meaning it only depends on the first $n$ batches $(I_0,\ldots,I_{n-1})$. 
The \textit{tower property} of conditional expectation, 
\[
\E[\phi] = \E[\E[\phi \vert \F_{n}] ] \text{ for any random variable } \phi \text{ on } (\Omega,\F,\P),
\]
will be frequently used.
To simplify the notation, for each index set $I=\lbrace i_{1},\ldots, i_{b} \rbrace \subset \lbrace 1,\ldots,p \rbrace$ of size $b$, we set
\[
Y_{I}:=Y_{i_{1}} \times \cdots \times Y_{i_{b}}, \quad y_{I}^{\delta} := \left( \yd_{i_{1}},\ldots,\yd_{i_{b}} \right), \quad \delta_{I}^{2} = \delta^{2}_{i_{1}} + \cdots + \delta^{2}_{i_{b}}
\]
and define $F_{I}:X \rightarrow Y_{I}$ by
\[
F_{I}(x):=\left( F_{i_{1}}(x),\ldots,F_{i_{b}}(x) \right), \quad \forall x \in X.
\]
Let $F^{\prime}_{I}$ denote the Fr\'{e}chet derivative of $F_{I}$. Then the method \eqref{xn} can be rewritten as
\[
\xd_{n+1} = \xd_{n} - \omega_n\frac{p}{b} F^{\prime}_{I_{n}}\left( \xd_{n} \right)^{\ast}\left(F_{I_{n}}(\xd_{n}) - \yd_{I_{n}} \right)
\]
with step control
\begin{equation*}
	\omega_n = \begin{cases}
		1 \quad \norm{F_{I_{n} }(\xd_{n}) - \yd_{I_{n}}}\geq\tau\delta_{I_{n}},\\
		0 \quad \text{else},
	\end{cases}
\end{equation*}
and the batched discrepancy principle simplifies to
\begin{equation*}
		n_\delta=\underset{n\geq L_\delta-1}{\arg\min}\Bigg\{ \sum_{j=0}^{L_\delta-1}\norm{F_{I_{n-j}} (\xd_{n-j}) - \yd_{I_{n-j}}}^2-\tau^2\delta_{I_{n-j}}^2\leq0\Bigg\}.
\end{equation*}
Note that $x_n^\delta$ is $\mathcal F_n$ measurable, but in general $F_{I_n}$ and $y_{I_n}$ are not, since the filtration $\mathcal F_n$ excludes the $(n+1)$-th batch $I_n$. This fact is used to calculate the expected values of the batch residual and the batch noise, which will be heavily used later on. 
\begin{lemma}\label{lem:Expec}
For an $\mathcal F_n$-measurable $x\in X$, we have
\begin{equation*}
\E\left[ \norm{\yd_{I_n} - F_{I_n}(x)}^{2} \right] = \frac{b}{p}\E\left[ \norm{\yd - F(x)}^{2}\right].
\end{equation*}
In particular, it holds
\begin{equation*}
\E\left[ \delta_{I_{n}}^{2} \right]  = \frac{b}{p}\delta^{2}.
\end{equation*}

\end{lemma}
\begin{proof} We follow the proof of \cite[Lemma 3.1]{Jin2023b}.
   Due to the $\mathcal F_n$-measurability of $x$, the random variable $\norm{\yd_{I_n} - F_{I_n}(x)}^{2}$ only depends on the batches $I_0,\ldots,I_n$.
    This implies that the conditional expectation $\Ex{\cdot \vert \mathcal F_n}$, which averages over all possible $I_l$ with $l\geq n$, corresponds only to the sum over all possible selections of $I_n$. Since $\mathbb P (I_n)=1/\binom p b$, it follows   
    \begin{align*}
\E\left[ \norm{\yd_{I_n} - F_{I_n}(x)}^{2}\Big\vert \F_{n} \right]
&= \frac{1}{\binom{p}{b}
} \sum_{I:\vert I \vert = b} \norm{\yd_{I} - F_{I}(x)}^{2}
\\&=\frac{1}{\binom{p}{b}
} \sum_{I:\vert I \vert = b} \sum_{i\in I} \norm{\yd_{i} - F_{i}(x)}^{2} \notag\\
&=\frac{1}{\binom{p}{b}} \sum_{i=1}^{p}  \sum_{I:\vert I \vert = b \, \& \, i \in I }\norm{\yd_{i} - F_{i}(x)}^{2}\\
&= \frac{\binom{p-1}{b-1}}{\binom{p}{b}
} \sum_{i=1}^{p} \norm{\yd_{i} - F_{i}(x)}^{2}
= \frac{b}{p}\norm{\yd - F(x)}^{2}\notag.
\end{align*}
In the third equation we swap the addends and the fourth equation follows from a combinatorial argument: in creating a subset $I$ containing a fixed $i$, there are $\binom{p-1}{b-1}$ ways to picks the remaining $b-1$ entries from the remaining $p-1$  equations. Applying the expectation on both sides and using the tower property leads us to the first result.\\
The second result follows analogously by replacing $F_{I_n}(x)$ with $y_{I_n}$.
\end{proof}
\begin{remark}
    Applying Lemma \ref{lem:Expec} to single-valued batches $I_n$ yields $\Ex{\delta_{I_n}^2}= \frac 1 p\delta^2$ and hence for $b$ independent single-valued batches it holds $\Ex{\delta_{I_n^{(1)}}^2+ \ldots+ \delta_{I_n^{(b)}}^2}=\frac b p \delta^2$. This argumentation implies that the result of Lemma \ref{lem:Expec} and therefore the rest of the paper can be generalized to the setting of choosing the batches with replacement. 
\end{remark}

To prove convergence we use the following standard Assumptions.
\begin{assumption}\label{A1}
Let the following conditions hold:
\begin{enumerate}[label=(\roman*)]
\item There is a $\rho>0$ such that problem \eqref{eq1} has a solution in $\mathcal{B}_{\rho}(x_{0})$ with $\mathcal B_{2\rho}(x_0)\subseteq\mathcal D(F)$.
\item Each operator $F_{i}:\mathcal{D}(F)\subset X \rightarrow Y_{i}$ has a continuous Fr\'{e}chet derivative $F_i^{\prime}$ on $\mathcal B_{2\rho}(x_0)\subseteq X$.
\item Each $F_{i}$ is properly scaled, so that $\Vert F_{i}^{\prime}(x) \Vert \leq 1/\sqrt{p} $ for all $x\in \mathcal B_{2 \rho}({x_0})$.
\item There exists an $\eta\in [0,1/2)$, such that the tangential cone condition
\[
\Vert F_{i}(x) - F_{i}(\bar{x}) - F^{\prime}_{i}(\bar{x})(x - \bar{x}) \Vert \leq \eta \Vert F_{i}(x) - F_{i}(\bar{x}) \Vert 
\]
holds for all $x,\bar{x} \in \mathcal{B}_{2\rho}(x_{0}) \subseteq \mathcal{D}(F) $ and for all $i=1,\ldots,p$.
\end{enumerate}
\end{assumption}
Assumption \ref{A1}(i) and (iv) guarantees the existence of the unique $x_0$-minimum norm solution $\xdag\in\mathcal{B}_{\rho}(x_{0})$, which is characterized by $\xdag-x_0\in \mathcal N(F^\prime(\xdag))^\perp$ \cite[Proposition 2.1]{Kaltenbacher2008}. Here $\mathcal N$ denotes the null space of an operator. Moreover, Assumption \ref{A1}(iv) is a rather strong condition since the nonlinearity condition is imposed to every function $F_{i}$ instead to the entire system via $F$. Due to the definition of norms on product spaces, it implies that
\begin{equation}\label{eq:resFn}
\Vert F_{I}(x) - F_{I}(\bar{x}) - F^{\prime}_{I}(\bar{x})(x - \bar{x}) \Vert \leq \eta \Vert F_{I}(x) - F_{I}(\bar{x}) \Vert 
\end{equation}
for all $x,\bar{x} \in \mathcal{B}_{2\rho}(x_{0}) \subseteq \mathcal{D}(F)$ and for every subset $I \subseteq \lbrace 1, \ldots, p\rbrace$. Accordingly, with $I = \lbrace 1, \ldots, p\rbrace$, the nonlinearity condition also holds for the overall operator $F:X\rightarrow Y$, i.e.
\begin{equation}\label{eq:resF}
\Vert F(x) - F(\bar{x}) - F^{\prime}(\bar{x})(x - \bar{x}) \Vert \leq \eta \Vert F(x) - F(\bar{x}) \Vert.
\end{equation}
Obviously, the reverse of the above implications is not true.\\
The following bound will be useful.
\begin{lemma}\label{lemma:Fbound}
Let Assumption \ref{A1}(iv) hold and $I\subseteq \{1,\ldots,p\}$. Then 
\begin{equation}\label{eq:Fbound}
\frac{1}{1+\eta}\norm{F_I^{\prime}(x)(x - \bar{x})} \leq\norm{F_I(x) - F_I(\bar{x})} \leq \frac{1}{1-\eta}\norm{F_I^{\prime}(x)(x - \bar{x})}
\end{equation}
for all $x,\bar x \in\mathcal{B}_{2\rho}(x_{0}) $.
\end{lemma}

\begin{proof}
Using equation \eqref{eq:resF} and the triangle inequality $\Vert a \Vert - \Vert b \Vert \leq \Vert a-b \Vert$ for $a,b \in Y$, we get
\begin{align*}
\eta \Vert F_I(x) - F_I(\bar{x}) \Vert&\geq\Vert F_I(x) - F_I(\bar{x}) - F_I^{\prime}(\bar{x})(x - \bar{x}) \Vert  \\&= \Vert F_I^{\prime}(\bar{x})(x - \bar{x}) - (F_I(x) - F_I(\bar{x}))   \Vert \\
&\geq \Vert F_I^{\prime}(\bar{x})(x - \bar{x}) \Vert - \Vert F_I(x) - F_I(\bar{x})\Vert,
\end{align*}
and hence,
\begin{align*}
\Vert F_I^{\prime}(\bar{x})(x - \bar{x}) \Vert &\leq (1+\eta) \Vert F_I(x) - F_I(\bar{x})\Vert\\
\frac{1}{1+\eta} \Vert F_I^{\prime}(\bar{x})(x - \bar{x}) \Vert &\leq \Vert F_I(x) - F_I(\bar{x})\Vert,
\end{align*}
which is the first inequality of \eqref{eq:Fbound}. Similarly it follows
\begin{align*}
 \Vert F_I(x) - F_I(\bar{x} )\Vert - \Vert F_I^{\prime}(x)(x- \bar{x}) \Vert &\leq \eta \Vert F_I(x) - F_I(\bar{x} )\Vert \\
 (1-\eta) \Vert F_I(x) - F_I(\bar{x}) \Vert &\leq \Vert F_I^{\prime}(x)(x- \bar{x}) \Vert \\
 \Vert F_I(x) - F_I(\bar{x}) \Vert &\leq \frac{1}{1-\eta} \Vert F_I^{\prime}(x)(x- \bar{x}) \Vert, \hspace{5em}
\end{align*}
which leads to the second inequality of \eqref{eq:Fbound}.
\end{proof}
The following result shows an error bound similar to that of Landweber iteration for nonlinear inverse problems derived in \cite[Proposition 2.2]{Hanke1995}. 
\begin{lemma}\label{lemma:err_ineq}
	Let Assumption \ref{A1} hold and $x_\ast\in \mathcal B _{2 \rho}(x_0)$ be an arbitrary solution of equation \eqref{eq1}. Consider the iteration $\{\xd_n\}\subseteq X$ generated by \eqref{xn} with a fixed selection of batches $(I_0,I_1,\ldots)\in \Omega$. Then for all iterates that satisfy $\xd_n\in \mathcal B_{2\rho}(x_0)$, it holds 
	\begin{multline*}
		\Vert\xd_{n+1} - x_{\ast}\Vert^{2} - \norm{\xd_{n} - x_{\ast}}^{2} \leq\omega_n\frac{p}{b}\norm{\yd_{I_{n}} - F_{I_{n}}(\xd_{n}) } \Big( 2(1 + \eta)\delta_{I_n} - (1-2\eta)\norm{\yd_{I_{n}} - F_{I_{n}}(\xd_{n}) }\Big).
	\end{multline*}
\end{lemma}
\begin{proof}
	We follow the same computations for the nonlinear Landweber iteration found in the classical paper of \cite{Hanke1995}, which also can also be read in \cite{Kaltenbacher2008}. We start with computing
	\begin{multline*}
		\norm{\xd_{n+1} - x_{\ast}}^{2} - \norm{\xd_{n} - x_{\ast}}^{2} = \norm{\xd_{n+1}}^2 - \norm{\xd_n}^2- 2 \langle \xd_{n+1},x_{\ast}\rangle + 2 \langle \xd_n, x_{\ast}\rangle\\
		= \norm{\xd_{n+1} - \xd_n}^2 + 2 \langle \xd_{n+1}, \xd_n\rangle - 2\norm{\xd_n}^2 - 2 \langle \xd_{n+1},x_{\ast}\rangle + 2 \langle \xd_n, x_{\ast}\rangle	\\
		= \norm{\xd_{n+1} - \xd_n}^2 +2 \langle \xd_{n+1} - \xd_n, \xd_n - x_{\ast}\rangle.
	\end{multline*}
	Inserting the definition of $\xd_{n+1}$ and using $\omega_n^2=\omega_n$, we obtain
	\begin{multline*}
		\norm{\xd_{n+1} - x_{\ast}}^{2} - \norm{\xd_{n} - x_{\ast}}^{2} \\= \omega_n\left(\frac{p}{b}\right)^{2} \norm{F^{\prime}_{I_{n}}(\xd_{n})^{\ast}(\yd_{I_{n}} - F_{I_{n}}(\xd_{n}))}^{2} + 2\omega_n\left(\frac{p}{b}\right)\langle \yd_{I_{n}} - F_{I_{n}}(\xd_{n}), F^{\prime}_{I_{n}}(\xd_{n})(\xd_{n} - x_{\ast}) \rangle.
	\end{multline*}
	To simplify the last expression and obtain a more useful estimate, we use the decomposition
	\begin{multline*}
		F^{\prime}_{I_{n}}(\xd_{n})(\xd_{n} - x_{\ast}) = (\yd_{I_{n}} -  y_{I_{n}}) + \left(F_{I_{n}}(x_{\ast}) - F_{I_{n}}(\xd_{n}) - F^{\prime}_{I_{n}}(\xd_{n})(x_{\ast} - \xd_{n}) \right) - (\yd_{I_{n}} - F_{I_{n}}(\xd_{n})).
	\end{multline*}
	Furthermore, by Assumption \ref{A1}(iii), it follows that
	\begin{equation*}
		\norm{F^{\prime}_{I_{n}}(\xd_{n})^{\ast}}^{2} = \norm{F^{\prime}_{I_{n}}(\xd_{n})}^{2} = \sum_{i \in I_n} \norm{F^{\prime}_{i}(\xd_{n})}^{2} \leq \frac{b}{p}.
	\end{equation*}
	Inserting this, we find 
	\begin{align*}
		\norm{\xd_{n+1} - x_{\ast}}^{2} &- \norm{\xd_{n} - x_{\ast}}^{2} \\
\leq&\omega_n \left(\frac{p}{b}\right) \norm{\yd_{I_{n}} - F_{I_{n}}(\xd_{n})}^{2} + 2\omega_n \left(\frac{p}{b}\right)\langle \yd_{I_{n}} - F_{I_{n}}(\xd_{n}), \yd_{I_{n}} -  y_{I_{n}} \rangle\\
		&+2\omega_n \left(\frac{p}{b}\right)\langle \yd_{I_{n}} - F_{I_{n}}(\xd_{n}) , F_{I_{n}}(x_{\ast}) - F_{I_{n}}(\xd_{n}) - F^{\prime}_{I_{n}}(\xd_{n})(x_{\ast} - \xd_{n})\rangle\\
		&-2\omega_n \left(\frac{p}{b}\right)\langle \yd_{I_{n}} - F_{I_{n}}(\xd_{n}), \yd_{I_{n}} - F_{I_{n}}(\xd_{n})  \rangle \\
		=&\omega_n\frac{p}{b}\Big(   2\langle \yd_{I_{n}} - F_{I_{n}}(\xd_{n}), \yd_{I_{n}} -  y_{I_{n}} \rangle\\
		&+2\langle \yd_{I_{n}} - F_{I_{n}}(\xd_{n}) , F_{I_{n}}(x_{\ast}) - F_{I_{n}}(\xd_{n}) - F^{\prime}_{I_{n}}(\xd_{n})(x_{\ast} - \xd_{n})\rangle- \norm{\yd_{I_{n}} - F_{I_{n}}(\xd_{n})}^{2} \Big).
	\end{align*}
	By the Cauchy-Schwarz inequality, we have
	\[
	2\langle \yd_{I_{n}} - F_{I_{n}}(\xd_{n}), \yd_{I_{n}} -  y_{I_{n}} \rangle
	\leq 2\norm{\yd_{I_{n}} - F_{I_{n}}(\xd_{n})} \norm{\yd_{I_{n}} -  y_{I_{n}}}
	\]
	Also, again by Cauchy-Schwarz, \eqref{eq:resFn} due to $\xd_n\in \mathcal B_{2 \rho}(x_0)$ and Assumption \ref{A1}(iv), we find
	\begin{align*}
		2\langle \yd_{I_{n}} - F_{I_{n}}(\xd_{n}) &,  F_{I_{n}}(x_{\ast}) - F_{I_{n}}(\xd_{n}) - F^{\prime}_{I_{n}}(\xd_{n})(x_{\ast} - \xd_{n})\rangle\\
		\leq& 2\norm{\yd_{I_{n}} - F_{I_{n}}(\xd_{n})}\norm{F_{I_{n}}(x_{\ast}) - F_{I_{n}}(\xd_{n}) - F^{\prime}_{I_{n}}(\xd_{n})(x_{\ast} - \xd_{n})}\\
		\leq& 2\eta \norm{\yd_{I_{n}} - F_{I_{n}}(\xd_{n})}\norm{F_{I_{n}}(x_{\ast}) - F_{I_{n}}(\xd_{n})}\\
		\leq&2\eta\norm{\yd_{I_{n}} - F_{I_{n}}(\xd_{n})} \norm{\yd_{I_{n}}- y_{I_{n}}} + 2\eta \norm{\yd_{I_{n}} - F_{I_{n}}(\xd_{n})}^{2}.
	\end{align*}
	This leads us to
	\begin{align*}
		\Vert\xd_{n+1} &- x_{\ast}\Vert^{2} - \norm{\xd_{n} - x_{\ast}}^{2} \notag \\
		\leq&\omega_n \left(\frac{p}{b}\right)\Big(  2\norm{\yd_{I_{n}} - F_{I_{n}}(\xd_{n})} \norm{\yd_{I_{n}} -  y_{I_{n}}} + 2\eta\norm{\yd_{I_{n}} - F_{I_{n}}(\xd_{n})} \norm{\yd_{I_{n}}- y_{I_{n}}} \\
		&+ 2\eta \norm{\yd_{I_{n}} - F_{I_{n}}(\xd_{n})}^{2} - \norm{\yd_{I_{n}} - F_{I_{n}}(\xd_{n}) }^{2}\Big) \notag \\
		=&\omega_n \left(\frac{p}{b}\right)\Big( 2(1 + \eta)\underbrace{\norm{\yd_{I_{n}}- y_{I_{n}}}}_{=\delta_{I_n}}\norm{\yd_{I_{n}} - F_{I_{n}}(\xd_{n}) } - (1-2\eta)\norm{\yd_{I_{n}} - F_{I_{n}}(\xd_{n}) }^{2}\Big),
	\end{align*}
	which is the desired result.
\end{proof}

\subsection{Convergence for exact data}
Next, we show the convergence of method \eqref{xn} with exact data. In this instance, we remove the superscript $\delta$ to denote exact data and iterates. Thus $\lbrace x_n \rbrace$ denotes the corresponding iterative sequence. We show that $\lbrace x_n \rbrace$ is a Cauchy sequence in $L^{2}(\Omega,X)$.
\begin{theorem}\label{thm:noisefree}
	Let Assumption \ref{A1} hold, $y\in Y$ be the exact data and $\{x_n\}$ the sequence of random variables on the probability space \eqref{eq:prob_space} generated by \eqref{xn} with exact data $y=\yd$. Then there exists a random variable $x^\ast \in L^2(\Omega,X)$, such that ${x^\ast}\in \mathcal B_{2 \rho} (x_0)$ solves problem \eqref{eq1} almost everywhere and $\{x_n\}$ converges in expectation to ${x^\ast}$, i.e.
	\[
	\E\left[ \norm{x_{n} - {x^\ast}} \right] \to 0 \quad \text{ as }n \to \infty.
	\]
	If $\mathcal{N}(F'(\xdag)) \subset \mathcal{N}(F'(x))$ for all $x \in \mathcal{B}_{\rho}(\xdag)$, then the sequence $\lbrace x_{n} \rbrace$ converges in expectation to the unique $x_{0}$-minimum-norm solution $\xdag$ of problem \eqref{eq1}:
	\[
	\E\left[ \norm{x_{n} - \xdag} \right] \to 0\quad \text{ as }n \to \infty.
	\]
\end{theorem}
\begin{remark}
Note that problem \eqref{eq1} can have multiple solutions and the random variable $x^\ast$ of Theorem \ref{thm:noisefree} is generally not constant. I.e. for the same starting value $x_0$ the method \eqref{xn} can converge to different solutions depending on the chosen batches within the iteration. This phenomenon only appears in the non-linear case, since for a linear operator $F$ it holds $\mathcal{N}(F'(\xdag)) = \mathcal{N}(F'(x))$ for all $x \in \mathcal{B}_{\rho}(\xdag)$ and hence the method converges to the unique $x_0$-minimum-norm solution.
\end{remark}
\begin{proof}[Proof of Theorem \ref{thm:noisefree}]
Let $x_\ast\in \mathcal B_{\rho}(x_0)$ be an arbitrary solution of problem \eqref{eq1} according to Assumption \ref{A1} (i). First we prove point-wise monotonicity of $\{\norm{x_n -x_\ast}\}$. For that, we fix an infinite selection of batches $\Lambda=(I_0,I_1,\ldots)\in \Omega$ and denote the corresponding iterates as $\hat x_n:=x_n(\Lambda)\in X$. Consequently, Lemma \ref{lemma:err_ineq} implies for all $\hat x_n\in \mathcal B_{2 \rho}(x_0)$ that
	\begin{equation}
    \label{eq:monotony_exact_data}
	\norm{\hat x_{n+1} -  x_{\ast}}^{2} - \norm{\hat x_{n} -  x_{\ast}}^{2}
	\leq-(1-2\eta)\frac{p}{b} \norm{y_{I_{n}} - F_{I_{n}}(\hat x_{n})}^{2}\leq 0.
	\end{equation}
    Here we used that $\delta_{I_n}=0$ and $\omega_n=1$ in the exact data case.
    Inductively, we get $\norm{\hat x_{n} - x_{\ast}}^{2}\leq\norm{x_{0} - x_{\ast}}^{2}<\rho^2$ and therefore $\hat x_n\in \mathcal{B}_\rho(x_\ast)\subseteq\mathcal B_{2 \rho } (x_0)$ for all $n\in \N$. Hence the inequality \eqref{eq:monotony_exact_data} holds for all $n\in \mathbb N_0$. And since $\Lambda \in \Omega$ was chosen arbitrarily the inequality \eqref{eq:monotony_exact_data} holds point-wise for $x_n$ and by taking expectations, we obtain
	\begin{align}
		\Ex{\norm{x_{n+1} - x_{\ast}}^{2}} - \Ex{\norm{x_{n} - x_{\ast}}^{2}} &\leq -(1-2\eta)\frac{p}{b} \Ex{\norm{y_{I_{n}} - F_{I_{n}}(x_{n})}^{2}} \notag \\
		&= -c\Ex{\norm{y - F(x_{n})}^{2}}, \label{eq:mono}
	\end{align}
    where we used the $\mathcal F_n$-measurability of $x_n$ and Lemma \ref{lem:Expec}, setting $c: = (1-2\eta)>0$.
Thus, the sequence $\lbrace \Exx{\norm{x_{n} - x_{\ast}}^{2}} \rbrace$ is monotonically decreasing and hence converges to some number $\alpha \geq 0$. Consequently, it also follows from \eqref{eq:mono} that
\begin{equation}\label{eq:Theta}
\lim_{n\to \infty}\Ex{\norm{y - F(x_{n})}^{2}}=0.
\end{equation}
In the following part we use the same method done in \cite{Hanke1995}. At this point, we need to show that the sequence $\lbrace x_{n} \rbrace$ is a Cauchy sequence in $L^{2}(\Omega,X)$. First note that $x_n\in L^2(\Omega,X)$ due to the above calculations. Now for every $j \geq n$, we choose $l\in \N$ with $j \geq l \geq n$ such that
\begin{equation}\label{eq:Ey}
\Ex{\norm{y - F(x_{l})}^{2}} \leq \Ex{\norm{y - F(x_{i})}^{2}}  \quad \text{ for all } n \leq i \leq j.
\end{equation}
Let $e_{n}:=x_{n} - x_{\ast}$. Using the fact that $\E\left[\Vert \cdot\Vert^{2}\right]^{1/2}$ is a norm \cite{Jin2020}, we get
\begin{equation}\label{eq:cauchyerr}
\Ex{\norm{e_{j} - e_{n}}^{2}}^{1/2} \leq \Ex{\norm{e_{j} - e_{l}}^{2}}^{1/2} + \Ex{\norm{e_{l} - e_{n}}^{2}}^{1/2}.
\end{equation}
Also in terms of expectation, based on the polarization identity, we have the identities
\begin{equation}\label{eq:errexpand}
\begin{aligned}
\Ex{\norm{e_{j} - e_{l}}^{2}}  &=  2\Ex{\left\langle e_{l} - e_{j} , e_{l} \right\rangle} + \Ex{\norm{e_{j}}^{2}}  - \Ex{\norm{e_{l}}^{2}},\\
\Ex{\norm{e_{l} - e_{n}}^{2}}  &=  2\Ex{\left\langle e_{l} - e_{n} , e_{l} \right\rangle} + \Ex{\norm{e_{n}}^{2}}  - \Ex{\norm{e_{l}}^{2}}.
\end{aligned}
\end{equation}
The last two terms in the right of both equations in \eqref{eq:errexpand} both tend to $\alpha-\alpha = 0$. We hence need to show that
\begin{equation}\label{eq6}
\Ex{\left\langle e_{l} - e_{n} , e_{l} \right\rangle} \to 0 \text{ as }n \to \infty.
\end{equation}
Now, by \eqref{xn:iteration},
\begin{align*}
\left\langle e_{l} - e_{n} , e_{l} \right\rangle
&=\sum_{k=n}^{l-1}
\left\langle e_{k+1}-e_{k}, e_{l} \right\rangle=-\sum_{k=n}^{l-1} \langle F_{I_{k}}(x_{k}) - y_{I_{k}}, F^{\prime}_{I_{k}}(x_{k})(x_{l} -  x_{\ast}) \rangle.
\end{align*}
By applying the triangle inequality and Cauchy-Schwarz inequality, we further get
\begin{align*}
\abv{\Ex{\left\langle e_{l} - e_{n} , e_{l} \right\rangle }}
\leq& \sum_{k=n}^{l-1}
\abv{\Ex{\langle F_{I_{k}}(x_{k}) - y_{I_{k}}, F^{\prime}_{I_{k}}(x_{k})(x_{l} - x_{k} + x_{k} - x_{\ast}) \rangle}}\\
=&\sum_{k=n}^{l-1}
\abv{\Ex{\langle F_{I_{k}}(x_{k}) - y_{I_{k}}, F^{\prime}_{I_{k}}(x_{k})(x_{l} - x_{k} ) \rangle}} \\
&+ \sum_{k=n}^{l-1}
\abv{\Ex{\langle F_{I_{k}}(x_{k}) - y_{I_{k}}, F^{\prime}_{I_{k}}(x_{k})( x_{k} - x_{\ast}) \rangle}} \\
\leq& \sum_{k=n}^{l-1}\Ex{\norm{F_{I_{k}}(x_{k}) - y_{I_{k}}}^{2}}^{1/2}\Ex{\norm{F^{\prime}_{I_{k}}(x_{k})(x_{l} - x_{k} ) }^{2}}^{1/2} \\
&+\sum_{k=n}^{l-1}\Ex{\norm{F_{I_{k}}(x_{k}) - y_{I_{k}}}^{2}}^{1/2}\Ex{\norm{F^{\prime}_{I_{k}}(x_{k})(x_{k} - x_{\ast} ) }^{2}}^{1/2}\\
&=:E_{1} + E_{2}.
\end{align*}
Applying Lemma \ref{lemma:Fbound} and the triangle inequality, we obtain
\begin{equation*}
\begin{aligned}
E_{1}
\leq& (1+\eta) \sum_{k=n}^{l-1}\Ex{\norm{F_{I_{k}}(x_{k}) - y_{I_{k}}}^{2}}^{1/2}\Ex{\norm{F_{I_{k}}(x_{l}) - F_{I_{k}}(x_{k})}^{2}}^{1/2}\\
\leq& (1+\eta) \sum_{k=n}^{l-1}\Ex{\norm{F_{I_{k}}(x_{k}) - y_{I_{k}}}^{2}}^{1/2}\Big(\Ex{\norm{F_{I_{k}}(x_{l}) - F_{I_{k}}(x_{\ast})}^{2}}^{1/2} \\
&\qquad\qquad + \Ex{\norm{F_{I_{k}}(x_{\ast}) - F_{I_{k}}(x_{k})}^{2}}^{1/2}\Big).
\end{aligned}
\end{equation*}
Since $x_{k}$ is $\F_{k}$-measurable, Lemma \ref{lem:Expec} implies
\[
\Ex{\norm{F_{I_{k}}(x_{k}) - y_{I_{k}}}^{2}} =  \frac{b}{p}\Ex{\norm{F(x_{k}) - y}^{2}}.
\]
Note that $x_{l}$ is not $\F_{k}$-measurable, since $k<l$. But observe that $\norm{F_{I_{k}}(x_{l}) - F_{I_{k}}(x_{\ast})} \leq  \norm{F(x_{l}) - F(x_{\ast})}$. By the choice of $l$ in \eqref{eq:Ey}, it holds that
\[
\Ex{\norm{F_{I_{k}}(x_{l}) - F_{I_{k}}(x_{\ast})}^{2}}^{1/2} \leq  \Ex{\norm{F(x_{l}) - F(x_{\ast})}^{2}}^{1/2} \leq \Ex{\norm{F(x_{k}) - F(x_{\ast})}^{2}}^{1/2}
\]
for $n \leq k < l$, noting that $y=F(x_{\ast})$. This brings us to
\[
E_{1} \leq c_{1}\sum_{k=n}^{l-1} \Ex{\norm{F(x_{k}) - F(x_{\ast})}^{2}},
\]
where $c_{1}:=(1+\eta)\sqrt{b/p}(1+\sqrt{b/p	})$. Similarly,
\begin{align*}
E_{2}
&\leq (1+\eta) \sum_{k=n}^{l-1}\Ex{\norm{F_{I_{k}}(x_{k}) - y_{I_{k}}}^{2}}^{1/2}\Ex{\norm{F_{I_{k}}(x_{k}) - F_{I_{k}}(x_{\ast})}^{2}}^{1/2}\\
&=(1+\eta) \sum_{k=n}^{l-1} \Ex{\norm{F_{I_{k}}(x_{k}) - y_{I_{k}}}^{2}} = c_{2}\sum_{k=n}^{l-1} \Ex{\norm{F(x_{k}) - y}^{2}},
\end{align*}
where $c_{2}:=(1+\eta)b/p$. Let $c_{3}:=c_{1}+c_{2}$. Combining the estimates for $E_{1}$ and $E_{2}$, and using \eqref{eq:mono}, we get
\begin{align*}
\abv{\Ex{\left\langle e_{l} - e_{n} , e_{l} \right\rangle }} &\leq c_{3}\sum_{k=n}^{l-1} \Ex{\norm{F(x_{k}) - y}^{2}} \leq c_{4}\sum_{k=n}^{l-1}\left( \Ex{\norm{e_{k}}^{2}} -  \Ex{\norm{e_{k+1}}^{2}} \right) \\
&= c_{4} \left(\Ex{\norm{e_{n}}^{2}} - \Ex{\norm{e_{l}}^{2}}\right),
\end{align*}
where $c_{4}:=\frac{c_{3}}{c}$. By a similar process, we can also deduce that
\[
\abv{\Ex{\left\langle e_{l} - e_{j} , e_{l} \right\rangle }} \leq c_{4} \left(\Ex{\norm{e_{l}}^{2}} - \Ex{\norm{e_{j}}^{2}}\right).
\]
Using the convergence of $\lbrace \Ex{\norm{e_n}^{2}} \rbrace$, it follows that the right hand sides of \eqref{eq:errexpand} converges to zero as $n \to \infty$. By \eqref{eq:cauchyerr}, it also follows that $\lbrace e_{n} \rbrace$ and, in turn, $\lbrace x_{n} \rbrace$ are Cauchy sequences in $L^{2}(\Omega,X)$. Thus $\brc{x_{n}}$ has a limit denoted by ${x^\ast}\in L^2(\Omega,X)$. Using \eqref{eq:Theta}, we get $\Exx{\norm{y - F(x^\ast)}^{2}}=0$ and, in turn, $x^\ast$ is a solution of problem \eqref{eq1} almost everywhere. The point-wise monotonicity of $\{\norm{x_n-x^\ast}\}$ and the convergence in expectation additionally imply $\norm{x_n-x^\ast}\overset{}{\to} 0$ almost everywhere as $n\to \infty$. Using this, we conclude from $x_{n}\in \mathcal B_{\rho}(x_ \ast)$, that ${x^\ast}\in \overline{\mathcal B_{\rho}(x_\ast)}\subseteq \mathcal B_{2\rho} (x_0)$ almost everywhere.\\
If, in addition, $\mathcal{N}(F'(\xdag)) \subset \mathcal{N}(F'(x))$ for all $x \in \mathcal{B}_{\rho}(\xdag)$, then one can simply apply the same argument as in \cite{Hanke1995}: Problem \eqref{eq1} has a unique solution $\xdag$ of minimal distance from initial guess $x_{0}$ which satisfies $\xdag - x_{0} \in \mathcal{N}(F^{\prime}(\xdag))^{\perp}$. We now again fix an infinite selection of batches $\Lambda=(I_0,I_1,\ldots)\in \Omega$ and denote the corresponding iterates as $\hat x_n:=x_n(\Lambda)\in X$. By definition of the iteration \eqref{xn:iteration}, it holds $$\hat x_{n+1} - \hat x_{n} \in R(F_{I_n}^\prime(\hat x_n)^\ast)\subset \mathcal{R}(F^{\prime}(\hat x_{n})^{\ast}) \subset \mathcal{N}(F^{\prime}(\hat x_n))^{\perp}\subset \mathcal{N}(F^{\prime}(\xdag))^{\perp}$$ for all $n\in \N_0$. Then it follows that $\hat x_{n}-x_{0} = (\hat x_{n} - \hat x_{n-1}) + \cdots + (\hat x_{1}-x_{0}) \in  \mathcal{N}(F^{\prime}(\xdag))^{\perp}$ for $n=0,1,\ldots$. Since $\Lambda \in \Omega$ and $n\in \N$ was arbitrary, this extends to the point-wise limit of $x_{n}$, i.e. ${x^\ast} - x_{0} \in \mathcal{N}(F^{\prime}(\xdag))^{\perp}$ almost everywhere, and therefore we get ${x^\ast} = \xdag$.
\end{proof}

\subsection{Convergence for noisy data}
\begin{lemma}\label{prop:DP}
Let Assumption \ref{A1} hold and $x_\ast\in \mathcal B_{2 \rho} (x_0)$ be any solution to problem \eqref{eq1}. Consider the iteration $\{\xd_n\}\subseteq X$ generated by \eqref{xn} with a fixed selection of batches $(I_0,I_1,\ldots)\in \Omega$ and
\begin{equation}\label{eq:tau}
\tau>\frac{2(1+\eta)}{1-2\eta}.
\end{equation}
Then $\xd_n\in \mathcal B_{2\rho}(x_0)$ for all $n\in \mathbb N_0$ and we get 
\begin{equation}\label{eq:monotonicity}
{\norm{\xd_{n+1} - x_{\ast}}}^2  - {\norm{\xd_{n} - x_{\ast}}}^2\leq-\omega_nc\norm{\yd_{I_{n}} - F_{I_{n}}(\xd_{n}) }^2,
\end{equation}
with a constant $c>0$. In particular, the sequence $\brc{\norm{\xd_{n} - x_{\ast}}}$ is monotonically decreasing.\end{lemma}
\begin{proof}
We follow the argumentation of \cite{Kaltenbacher2008}[Proposition 2.2].
Let us denote the constant $c_1\in \mathbb R$ as $c_{1}:=(1-2\eta) - 2(1+\eta)\tau^{-1}$. Note that $c_{1}>0$ by \eqref{eq:tau}. Suppose $\xd_n\in\mathcal B_{2\rho}(x_0)$ for a fixed $n\in \mathbb N_0$. Then Lemma \ref{lemma:err_ineq} implies
\begin{align}
\Vert\xd_{n+1} &- x_{\ast}\Vert^{2} - \norm{\xd_{n} - x_{\ast}}^{2} \notag\\
&\leq \omega_n\frac{p}{b}\norm{\yd_{I_{n}} - F_{I_{n}}(\xd_{n}) }\Big( 2(1 + \eta)\delta_{I_n} - (1-2\eta)\norm{\yd_{I_{n}} - F_{I_{n}}(\xd_{n}) }\Big)  \notag \\
&\leq -c_1\omega_n\frac{p}{b}\norm{\yd_{I_{n}} - F_{I_{n}}(\xd_{n}) }^{2}.\label{eq:errortau}
\end{align}
The last inequality follows directly by the definition of $\omega_n$: If $\omega_n=0$ both sides vanish and if $\omega_n=1$ it holds $\delta_{I_n}\leq \tau^{-1}\norm{\yd_{I_{n}} - F_{I_{n}}(\xd_{n}) } $.\\
In particular, for any solution $x_\ast\in \mathcal B_\rho (x_0) $ according to Assumption \ref{A1}(i), it follows inductively $\xd_{n}\in \mathcal B_\rho (x_\ast)\subseteq \mathcal B_{2 \rho}(x_0) $ for all $n\in \N_0$, since $\xd_{n-1}\in\mathcal B_\rho (x_\ast)$ implies $\xd_{n}\in\mathcal B_\rho (x_\ast)$ by $$\norm{\xd_{n} - x_{\ast}} \leq \norm{\xd_{n-1} - x_{\ast}} <\rho.$$  
Hence inequality \eqref{eq:errortau} holds for all iterates with $c=c_1 p/b$. 
\end{proof}
In particular, the above Lemma implies that $\xd_n\in L^2(\Omega,X)$ for the random variable $\xd_n$ generated by \eqref{xn} with $\tau$ as in \eqref{eq:tau}. In the next Lemma, we analyze the relationship between the noisy and the noise-free iteration.
\begin{lemma}\label{lem:stability_batchDP_stepsize}
Let Assumption \ref{A1} hold and $N\in \N_{0}$ be fixed. Consider the iteration $\{\xd_n\}\subseteq X$ defined by \eqref{xn} with an arbitrary selection of batches $(I_0,I_1,\ldots) \in \Omega$. Then it holds 
\begin{equation*}
    \norm{\xd_N - x_N} \to 0 \text{ as } \delta \to 0.
\end{equation*}
In particular, for the sequence of random variables $\{\xd_n\}\subseteq L^2(\Omega,X)$, it holds
\[
 \Ex{\norm{\xd_{N} - x_{N}}} \to 0 \text{ as } \delta \to 0.
\]
\end{lemma}
\begin{proof}
The proof follows the argumentation of \cite[Lemma 3.3]{Jin2023b}.
First note that the convergence of $\xd_N$ for a fixed selection of batches implies the convergence in expectation, since there are only finitely many sample paths for fixed $N\geq 0$.\\
We prove by mathematical induction over $N$. For $N=0$ the result is trivial. Assume the result holds for $\xd_N$ with a fixed $N\in \N$, we need to prove it holds for $\xd_{N+1}$. Now, by \eqref{xn:iteration}, and using $\omega_N=1$ in the exact data case, it holds
\begin{align*}
\xd_{N+1} - x_{N+1} 
=& (\xd_{N} - x_{N}) - \omega_N\frac{p}{b} F^{\prime}_{I_{N}}(\xd_{N})^{\ast}(F_{I_{N}}(\xd_{N}) - \yd_{I_{N}})+\frac{p}{b} F^{\prime}_{I_{N}}(\xd_{N})^{\ast}(F_{I_{N}}(x_{N}) - y_{I_{N}})\\
&-\frac{p}{b}F^{\prime}_{I_{N}}(\xd_{N})^{\ast}(F_{I_{N}}(x_{N}) - y_{I_{N}})+ \frac{p}{b}F^{\prime}_{I_{N}}(x_{N})^{\ast}(F_{I_{N}}(x_{N}) - y_{I_{N}})\\
=&(\xd_{N} - x_{N}) - \frac{p}{b}\left( F^{\prime}_{I_{N}}(\xd_{N})^{\ast} - F^{\prime}_{I_{N}}(x_{N})^{\ast}\right)\big(F_{I_{N}}(x_{N}) - y_{I_{N}} \big) \\
&+\frac{p}{b} F^{\prime}_{I_{N}}(\xd_{N})^{\ast}\big(\omega_N(F_{I_{N}}(\xd_{N}) - \yd_{I_{N}}) - (F_{I_{N}}(x_{N}) - y_{I_{N}})   \big).
\end{align*}
Then by triangle and Cauchy-Schwarz inequality,
\begin{align*}
\norm{\xd_{N+1} - x_{N+1}} \leq& \norm{\xd_{N} - x_{N}} + \frac{p}{b}\norm{F^{\prime}_{I_{N}}(\xd_{N})^{\ast} - F^{\prime}_{I_{N}}(x_{N})^{\ast} }\norm{F_{I_{N}}(x_{N}) - y_{I_{N}} } \\
&+\frac{p}{b}\norm{F^{\prime}_{I_{N}}(\xd_{N})^{\ast}}\norm{\omega_N(F_{I_{N}}(\xd_{N}) - \yd_{I_{N}}) - (F_{I_{N}}(x_{N}) - y_{I_{N}}) }.
\end{align*}
By induction hypothesis the first term converges to zero for $\delta\to 0$. Using the continuity of the operator $F_I^\prime$ the same holds for the second term. For the third term first note that with the same argumentation it follows $\limsup_{\delta \to 0} \norm{F^{\prime}_{I_{N}}(\xd_{N})^{\ast}}<\infty$. We now want to show that 
\begin{equation}
\label{eq:noise_residuum_error}
    \norm{\omega_n (F_{I_{N}}(\xd_{N}) - \yd_{I_{N}}) - (F_{I_{N}}(x_{N}) - y_{I_{N}})  }\to 0 \text{ as } \delta \to 0.
\end{equation}
For this, let $R(\delta):=F_{I_{N}}(\xd_{N}) - \yd_{I_{N}}$ and $R(0):=F_{I_{N}}(x_{N}) - y_{I_{N}}$ its exact data variant. This reduces \eqref{eq:noise_residuum_error} to $\lim_{\delta \to 0} \norm{\omega_N  R(\delta) - R(0) }=0$ with $\omega_N=1$ if $\norm{R(\delta)}\geq \tau \delta$ and $\omega_N=0$ else. Since $\lim _{\delta \to 0} y_{I_N}^\delta =y_{I_N} $, we get with the induction hypothesis and continuity of $F_{I_N}$ that $\lim _{\delta \to 0 } R(\delta) = R(0)$. For $R(0) =0$ this implies \eqref{eq:noise_residuum_error}. If $\norm{R(0)}>0$ then $\norm{R(\delta)} \geq \tau \delta$ and hence $\omega_N=1$ for sufficiently small $\delta$. Therefore \eqref{eq:noise_residuum_error} and the statement follows.
\end{proof}

In the following Theorem we prove a general result, which includes the convergence for both the a priori and a posteriori stopping.
\begin{theorem}
\label{thm:general_convergence}
Let Assumption \ref{A1} hold and choose $\tau$ according to \eqref{eq:tau}. Consider the batched Landweber-Kaczmarz iteration $\{\xd_n\}\subseteq L^2(\Omega,X)$ defined by \eqref{xn} and let $x^\ast \in L^2(\Omega,X)$ be the random variable from Theorem \ref{thm:noisefree}, which solves problem \eqref{eq1} almost everywhere. Then for all random variables $n_\delta\in \N$ on the probability space \eqref{eq:prob_space}, with $n_\delta \to \infty$ as $\delta \to 0$, it holds
    \begin{equation*}
        \Ex{\norm{\xd_{n_\delta}- {x^\ast} }}\to 0 \text{ as } \delta \to 0,
    \end{equation*}
    where $\E$ denotes the joint expectation.
\end{theorem}
\begin{proof}
First, note that the monotonicity of the error $\{\norm{x_n - {x^\ast}}\}$ (Lemma \ref{prop:DP}) and its convergence in $L^2(\Omega,X)$ (Theorem \ref{thm:noisefree}) imply the point-wise convergence $\norm{x_n - {x^\ast}}\to 0 $ almost everywhere as $n\to \infty$. 
Let now $\{y^{\delta_k}\}\subseteq Y$ be a sequence of noisy data with $\lim_{k \to \infty}\delta_k = 0$. We now prove point-wise convergence of the noisy iterates, i.e. $\normm{x^{\delta_k}_{n_{\delta_k}}-x^\ast}\overset{}{\to} 0$ almost everywhere as $k\to \infty$.  
 For this, we fix an infinite selection of batches $\Lambda:=(I_0,I_1,\ldots)\in \Omega$, such that $\hat x^\ast:= x^\ast(\Lambda)$ is a solution of \eqref{eq1} and for the iteration $\hat x_n:=x_n(\Lambda)$ it holds $\lim_{n\to \infty}\norm{\hat x_{n}-{\hat x^\ast}}=0$.\\
Setting $n_k:=n_{\delta_k}(\Lambda) > m$, Lemma \ref{prop:DP} implies that
	\begin{equation}
    \label{eq:E(n)}
	\normm{\hat x^{\delta_{k}}_{n_k } - {\hat x^\ast}} \leq \normm{\hat x^{\delta_{k}}_{m} - {\hat x^\ast}} \leq  \normm{\hat x^{\delta_{k}}_{m} - \hat x_{m}} + \normm{\hat x_{m} - {\hat x^\ast}}.
	\end{equation}
	Let $\varepsilon>0$. Since the second summand in \eqref{eq:E(n)} vanishes for large $m$, we can fix some  $m=m(\varepsilon)$ so that the term is less than $\varepsilon/2$. 
	Also, by Lemma \ref{lem:stability_batchDP_stepsize}, there is $k_0=k(\varepsilon,m)$ such that the first term in the right hand side of \eqref{eq:E(n)} is less than $\varepsilon/2$ for $k \geq k_0$. And due to $n_k \to  \infty$ there is a $k_1(m)$, such that $n_k >m$ for all $k\geq k_1$. Consequently for all $k\geq \max\{k_0,k_1\}$ we get $\normm{\hat x^{\delta_{k}}_{n_{k}} - \hat x^{\ast}}\leq \varepsilon$ and hence $\lim_{k\to \infty}\normm{\hat x^{\delta_{k}}_{n_k} - \hat x^{\ast}}=0$. \\
    Since the selection of the batches $\Lambda$ was arbitrary up to a null set, this implies the convergence in expectation using dominated convergence with majorant $\norm{x_0 - x^\ast}$, since $$\norm{\xd_{n_\delta} - x^\ast}\leq\norm{\xd_{0} - x^\ast}= \norm{x_0 - x^\ast}<\infty.$$
\end{proof}
\begin{corollary}[A priori convergence]
\label{cor:a_priori}
Let the assumptions of Theorem \ref{thm:general_convergence} hold. Then for all $n_\delta\in \mathbb N$ with $n_\delta\to \infty$ as $\delta \to 0$, it holds
        \begin{equation*}
        \Ex{\norm{\xd_{n_\delta}- {x^\ast} }}\to 0 \text{ as } \delta \to 0.
    \end{equation*}
\end{corollary}
\begin{remark}
\label{remark:why_early_stoping}
    Note that in Theorem \ref{thm:general_convergence} and Corallary \ref{cor:a_priori} no growth condition on $n_\delta$ is needed. This is a result of the step size $\omega_n$ in \eqref{xn:step_control}. Unlike the typical semiconvergent behaviour of the classical Landweber iteration, iteration \eqref{xn} yields unconditional monotonicity of $\norm{\xd_n - x_\ast}$. Instead of deteriorating the approximation for large $n$ the iteration remains constant. However, a larger stopping index $n_\delta$ corresponds to a larger computational effort. Hence even though the size of $n_\delta$ has theoretically no impact, early stopping is still computationally desirable. 
\end{remark}

\subsubsection{Batched Discrepancy Principle}
\label{section:batch_DP}
Let us now analyze the batched discrepancy principle \eqref{eq:summed_batch_DP}.
The main goal of the proposed criterion is to approximate the classical discrepancy principle with a computationally cheaper rule. For this one observes that  $\norm{F_{I_n }(\xd_n) - \yd_{I_n}}^2-\tau^2 \delta_{I_n}^2$ is an unbiased estimator of \linebreak $\frac bp \left(\norm{F(\xd_n) - \yd}^2- \tau ^2\delta^2\right)$ by Lemma \ref{lem:Expec}. We improve this estimate by aggregation and approximate 
\begin{equation*}
		 \frac bp \left(\norm{F(\xd_n) - \yd}^2- \tau ^2\delta^2\right) \approx \frac 1 {L_\delta}\sum_{j=0}^{L_\delta-1}\left( \norm{F_{I_{n-j}} (\xd_{n-j}) - \yd_{I_{n-j}}}^2-\tau^2\delta_{I_{n-j}}^2 \right), 
\end{equation*}
where the \textit{aggregation length} $L_\delta \in \mathbb N$ is fixed a priori. Here we ignore the errors due to evaluation at different iterates, since we expect them to be similar for  $L_\delta$ not too large. However, for smaller $L_\delta$ less equations are monitored and the probability of overlooking equations in an iteration increases, i.e. it is more likely that $i\notin \bigcup_{j=0}^{L_\delta -1} I_{n-j}$ for an equation $i\in\{1,\ldots,p\}$. This might result in too early stopping, leading to poor reconstructions and no convergence in expectation. \\ Overall, the intuition suggests that the best results can be obtained for a sufficiently large but not too large aggregation length $L_\delta$. We will see later that the parameter $L_\delta$ does indeed have a significant impact on the quality of the reconstructions.\\Before proving convergence, the next Lemma shows that the proposed rule is well-defined, i.e. the iteration stops. 
\begin{lemma}
\label{lemma:finite_termination}
    Let Assumption \ref{A1} hold, assume $\delta_i>0$ for all $i\in \{1,\ldots,p\}$, $L_\delta \in \mathbb N$ and choose $\tau>1$ according to \eqref{eq:tau}. Then iteration $\{\xd_n\}\subseteq X$ generated by \eqref{xn} with an arbitrary selection of batches $(I_0,I_1,\ldots)\in \Omega$ terminates under the batched discrepancy principle \eqref{eq:summed_batch_DP} in a finite number of iterations, i.e. $\nd<\infty$.
\end{lemma}
\begin{proof}
Let $x_\ast\in \mathcal B_\rho(x_0)$ be a solution of problem \eqref{eq1} according to Assumption \ref{A1} (i). With $c>0$ from Lemma \ref{prop:DP}, it follows for all $M\in \mathbb N$ that
\begin{align*}
c\sum_{n=0}^{M}\omega_n\norm{\yd_{I_{n}} - F_{I_{n}}(\xd_{n})}^2  &\leq \sum_{n=0}^{M}\Vert\xd_{n} - x_{\ast}\Vert^{2}  - \Vert \xd_{n+1} - x_{\ast}\Vert^{2}\\
&\leq \Vert\xd_{0} - x_{\ast}\Vert^{2} - \Vert\xd_{M+1} - x_{\ast}\Vert^{2}\\
&\leq \Vert x_{0} - x_{\ast}\Vert^{2}  < \infty.
\end{align*}
In particular for $M\to \infty$ this implies 
\begin{equation*}
    \omega_n\norm{\yd_{I_{n}} - F_{I_{n}}(\xd_{n})}^2\to 0 \text{ as } n \to \infty.
\end{equation*}
Since $\delta_i>0$ for all $i\in \{1,\ldots,p\}$, there exists an $N\in \mathbb N$ such that
\begin{equation*}
    \omega_n \norm{\yd_{I_{n}} - F_{I_{n}}(\xd_{n})} ^2\leq   \min_{i\in \{1,\ldots,p\}} \tau^2\delta_i^2 \leq \tau ^2\delta_{I_n}^2 \quad \text{for all } n>N.
\end{equation*}
By multiplying both sides with $\omega_n$ and adding up $L_\delta$ consecutive iterates, it follows
\begin{equation*}
    \sum_{l=0}^{L_\delta-1}\omega_{n-l}\norm{F_{I_{n-l}}(\xd_{n-l}) - \yd_{I_{n-l}} }^2 \leq \sum_{l=0}^{L_\delta-1} \omega_{n-l}\tau^2\delta_{I_{n-l}}^2\quad \text{for all } n>N-L_\delta -1.
\end{equation*}
Note that $\omega_{n-l}^2=\omega_{n-l}$, since it only takes values $0$ and $1$. Additionally we get
\begin{equation*}
    \norm{F_{I_{n-l} }(\xd_{n-l}) - \yd_{I_{n-l}}}^2-\tau^2 \delta_{I_{n-l}}^2\leq  \omega_n\left(\norm{F_{I_{n-l} }(\xd_{n-l}) - \yd_{I_{n-l}}}^2-\tau^2 \delta_{I_{n-l}}^2\right)
\end{equation*}
for all $n\geq l$, because per definition $\omega_{n-l}=0$ iff the right hand side is negative and $\omega_{n-l}=1$ else. Together this results in
\begin{equation*}
    \sum_{l=0}^{L_\delta-1} \left( \norm{F_{I_{n-l} }(\xd_{n-l}) - \yd_{I_{n-l}}}^2-\tau^2 \delta_{I_{n-l}}^2 \right) \leq \sum_{l=0}^{L_\delta-1} \omega_n\left(\norm{F_{I_{n-l} }(\xd_{n-l}) - \yd_{I_{n-l}}}^2-\tau^2 \delta_{I_{n-l}}^2\right)\leq 0
\end{equation*}
for all $n>N+{L_\delta}-1$. Hence $\nd\leq N+L_\delta<\infty$ which proves the statement.
\end{proof}
\begin{theorem}[A posteriori convergence]
\label{thm:batch_DP_conv}
Let the assumptions of Theorem \ref{thm:general_convergence} hold and consider $L_\delta\in \N$, such that $L_\delta \to \infty$ as $\delta \searrow 0$. Then for $n_\delta\in \N$ chosen according to \eqref{eq:summed_batch_DP}, it holds
        \begin{equation*}
        \Ex{\norm{\xd_{n_\delta}- {x^\ast} }}\to 0 \text{ as } \delta \searrow 0.
    \end{equation*}
\end{theorem}
\begin{proof}
    This is a direct consequence of Theorem \ref{thm:general_convergence}, since $\infty>n_\delta\geq L_\delta \to \infty$ as $\delta \searrow 0$.
\end{proof}
The condition $L_\delta \to \infty$ is necessary to ensure that the probability of overlooking equations vanishes. As in the previous results, there is no growth condition needed on $L_\delta$ due to the step control $\omega_n$ in the definition of the method \eqref{xn}. However, in the next section will see, that it is well advised to control the parameter $L_\delta$. 

\section{Numerical illustration}\label{Section:numerics}
In this section, the batched Kaczmarz-Landweber method is tested on a simulation of the electrical impedance tomography problem (EIT-problem). EIT is a type of medical imaging to determine the electrical conductivity inside a domain by boundary measurements of currents and voltages. The formulation in this paper is based on \cite{Borcea_2002}. Consider the conductivity equation
\begin{equation*}
    -\nabla \cdot( \sigma \nabla u) =0. 
\end{equation*}
In the physical context, $\sigma $ denotes the conductivity and $u$ the voltage on a domain $\domain \subseteq \mathbb R^d$.
The theoretical EIT-problem, also known as Calderón problem, is the parameter identification of the conductivity $\sigma$ from the Neumann-to-Dirichlet map (N2D-map).

The N2D-map $\Lambda_\sigma(g)$ is defined as boundary evaluation of the solution of the conductivity equation with Neumann boundary value $g$. I.e. 
$\Lambda_\sigma(g):= u_{\sigma,g}\vert_{\partial \domain}$ where $u_{\sigma,g}$ is the unique solution of
\begin{equation}
\label{equation:strong_conductivity}
\begin{cases}
    -\nabla \cdot(\sigma \nabla u) =0\quad &\text{in} \; \domain\\
    \sigma \frac {\partial u} {\partial n} = g \quad &\text{on} \; \partial \domain.
    \end{cases}
\end{equation}
Hence the forward operator of the EIT-problem $F$ maps the conductivity $\sigma$ to the N2D-map  $F: \sigma \to \Lambda _\sigma$.\\ 
While the just defined forward operator corresponds to the problem with all Neumann input functions, a single measurement is given by $F_f(\sigma):=F(\sigma)f$ for Neumann input data $f$. For multiple Neumann input data $f_1,\ldots,f_p$ the inverse problem corresponds to the prior introduced setting
\begin{equation*}
    F_{f_i}(\sigma) = y_i \quad i=1,\ldots,p,
\end{equation*}
with noisy data $y_i^\delta$ satisfying
\begin{equation*}
    \norm{y_i^\delta -y_i} \leq \delta_i.
\end{equation*}
To formally define the operator we introduce the spaces with vanishing boundary integral $\Hdmnd{1}\domain$ and $\LLdmnd{\partial \domain}$.
 Using this spaces, the weak formulation of \eqref{equation:strong_conductivity} is given by 
\begin{equation}
    \label{equation:weak_conductivity_equation}
    \int_\domain \sigma \nabla u \nabla v \mathrm d x = \int_{\partial \domain} g \mathrm{Tr} v \mathrm d S \quad \forall v\in \Hdmnd 1 \domain,
\end{equation}
where $\mathrm {Tr}: \Hs 1 \domain \to \Lp 2 {\partial \domain}$ denotes the trace operator.
We define the corresponding differential operator $L_\sigma: \Hdmnd 1 \domain \to \Hdmnd{-1}{\domain}$ as $L_\sigma(u)= \int_\domain \sigma \nabla u \nabla \cdot \mathrm d x$. If $\sigma \in \DefSigma$, with $$\DefSigma:= \{\sigma \in \Lp \infty \domain: \mathrm{ess} \inf \sigma>0 \}.$$ it can be shown, analogously to the Dirichlet setting in \cite[Section 6.2.1]{evans2010partial}, by using the Lax-Milgram theorem and the Poincaré inequality, that $L_\sigma$ is boundedly invertible.\\ 
Using the dual map of the trace operator $\mathrm{Tr} ^\prime : \Lp 2 {\partial \domain} \to \Hdmnd {-1} \domain$ the weak conductivity equation with Neumann boundary values  \eqref{equation:weak_conductivity_equation}  reduces to $L_\sigma u = \mathrm{Tr}^\prime g$ and hence has the unique solution $u_{\sigma,g}=L_\sigma^{-1}\mathrm{Tr}^\prime g$ for $\sigma \in \DefSigma$.\\
Therefore, the forward operator of the EIT-problem is formally defined as
\begin{equation*}
    F:\DefSigma\subseteq \Lp\infty\domain \to \mathcal L(\LLdmnd{\partial \domain},\LLdmnd{\partial \domain}) \quad \text{with} \quad    F(\sigma) = \mathrm {Tr} L_\sigma^{-1} \mathrm{Tr}^\prime .
\end{equation*}
Using basic calculus we can show that the operator is Frechét-differentiable with respect to the $\mathbf L^\infty$-norm and the derivative is given by 
\begin{equation*}
    F^\prime(\sigma)h= - \mathrm {Tr} L_\sigma ^{-1} L_h L_\sigma^{-1} \mathrm{Tr}^\prime. 
\end{equation*}
However, the considered batched Landweber-Kaczmarz method requires the operator to be defined on a Hilbert space and the Frechét-derivative does not extend to $\mathbf L^2$. Therefore we will work with Sobolev spaces $\Hs s \domain$ based on the Fourier transform with $s>0$. This spaces are Hilbert spaces and continuously embed into $\Lp \infty \domain$ for $s>d/2$ where $d$ is the dimension of the underlying domain $\domain$, see e.g. \cite[Theorem 1.66]{Bahouri_2011}. Let now $s>d/2$, then the restricted operator $F \vert _{\Hs s \domain}$ is Frechét-differentiable with respect to $\mathbf H^s$. Analogously the restricted forward operator with fixed Neumann data $F_f:\DefSigma\cap \Hs s \domain \subseteq \Hs s \domain \to \LLdmnd{\partial \domain}$, defined as $F_f(\sigma)h=F(\sigma)f$, is Frechét-differentiable with derivative
\begin{equation}\label{eq:Frechet_f}
    F_f^\prime(\sigma)h= - \mathrm {Tr} L_\sigma ^{-1} L_h L_\sigma^{-1} \mathrm{Tr}^\prime f.
\end{equation}
The derivative can be interpreted as boundary evaluation $F_f^\prime(\sigma)h = \mathrm{Tr} \tilde u$ of the unique solution $\tilde u\in \Hdmnd{1} \domain$ of the boundary value problem
\begin{equation*}
\begin{cases}
     -\nabla\cdot (\sigma  \nabla \tilde u) =  \nabla \cdot (h \nabla u_{\sigma,f}) \quad &\text{in} \; \domain\\
    \sigma \frac {\partial \tilde u} {\partial n} = -h\frac {\partial u_{\sigma,f}} {\partial n} \quad &\text{on} \; \partial \domain.
    \end{cases}
\end{equation*}
Here again $u_{\sigma,f} = L_\sigma^{-1}\mathrm{Tr}^\prime f\in \Hdmnd 1 \domain$ is the solution to the conductivity equation \eqref{equation:weak_conductivity_equation} with boundary value $f$.\\
Since the derivative $F^\prime_{f}$ is continuous, the operator $F_f$ fulfills Assumption \ref{A1} (ii). The boundedness of the Frechét-derivative (Assumption \ref{A1} (iii)) is ensured by scaling the operator $F_f$ properly. This means applying the method on adjusted data $\hat y^\delta:=(\rho_1y_{f_i}^\delta,\ldots, \rho_p y_{f_p}^\delta)$, using the scaled operator $\hat F:=(\rho_1 F_{f_1},\ldots,\rho_p F_{f_p})$ with 
\begin{equation}
\label{eq:EIT_scaling_factor}
\rho_i = \frac 1 {\sqrt p}\cdot\frac 1 {  \max\{c_F,\sup_{\sigma \in \DefSigma^s} ||F^\prime_{f_i}[\sigma]||_{\mathcal L(\Hs s \domain , \LLdmnd{\partial \domain})}\}}
\end{equation}
where $c_F>0$ is a positive constant. 
For the noise level $\hat \delta= (\hat \delta_1,\ldots,\hat \delta_p)$, this implies
$$
\norm{\hat y_i - \hat y_i^\delta}\leq \rho_i \delta_i =: \hat \delta_i \quad i=1,\ldots,p.
$$
The validity of Assumption \ref{A1} (iv), the tangential cone condition for EIT, is an open problem and only known in special cases (see \cite{Kindermann_2022} for a discussion).\\
To implement the Landweber-Kaczmarz iteration we use a discretization of the continuous Landweber step, i.e. evaluate the operator
\begin{align*}
    \sigmad_n\to \sigmad_n - \omega_n\frac p b \sum_{i\in I_n}  \hat F^\prime_{f_i}(\sigmad_n)^\ast ( \hat F_{f_i}(\sigmad_n) - \hat y_i^\delta).
\end{align*}
For this we calculate the adjoint $F_f^\prime(\sigma)^\ast$ as
\begin{equation*}
 F_f^\prime(\sigma)^\ast g = - \mathcal I_{\Hs s \domain}^{-1}\left(\langle\nabla u_{\sigma,f} \nabla u_{\sigma,g} , \cdot\rangle_{\Lp 2 \domain}\right).
\end{equation*}
Here again $u_{\sigma,h}=L_\sigma^{-1}\mathrm{Tr}^\prime h$ is the solution of the conductivity equation with Neumann boundary value $h\in \LLdmnd{\partial \domain}$ (cf. \eqref{equation:weak_conductivity_equation}) and $\mathcal I_{\Hs s \domain}: \Hs s \domain \to \Hs s \domain ^\ast$ denotes the Riesz-operator of $\Hs s \domain$. The implementation uses the finite element method with p1-triangles to solve the boundary value problem and the discrete Fourier transform to approximate the Riesz-operator of $\Hs s \domain$. Note that the discretization of the continuous adjoint $F^\prime_f(\sigma)^\ast$ is just an approximation of the actual discrete adjoint.\\

We compare the batched discrepancy principle, introduced in this paper, with the classical discrepancy principle \eqref{landwDP}.
To  simplify the notation of both termination criteria we define the  discrepancy ratios as
\begin{align*}
    \DiscRatio_{n}^{b,L} := \frac {\sum_{j=0}^{L-1} \normm{\hat F_{I_{n-j}}(\sigmad_{n-j})-\hat y^\delta_{I_{n-j}}}^2}{\sum_{j=0}^{L-1}\hat\delta_{I_{n-j}}^2}
    \quad\text{and}\quad
    \DiscRatio_{n}^{\text{full}} :=  \frac {\normm{\hat F(\sigmad_n) - \hat y^\delta}^2} {{\hat\delta}^2}.
\end{align*}
Using this notation, the iteration terminates at the first $n\in \mathbb N$ with $\DiscRatio_n< \tau ^2$ where $\DiscRatio_n\in \{\DiscRatio_n^{b,L},\DiscRatio_n^{\text{full}}\}$.
The final method is given in Algorithm \ref{kaczmarz_algo_eit}. Note that the scaling factor $\rho_i$ of equation \eqref{eq:EIT_scaling_factor} is simplified by using the operator norm of the current term $F_{f_i}^\prime(\sigmad_n)$ instead of the supremum over all feasible conductivities.
\begin{algorithm}
\caption{Landweber-Kaczmarz algorithm for EIT}
\label{kaczmarz_algo_eit}
\begin{algorithmic}[1]
    \State Choose a domain discretization (mesh and grid).
    \State Select boundary functions $\mathbf{f} = (f_1, \ldots, f_p)$ and pick an initial guess $\sigma_0 \in \DefSigma$.
    \State Choose $\tau > 1$ and $c_F>0$ appropriately.
    \State Select the discrepancy type $\DiscRatio_n \in \{\DiscRatio_n^{b,L}, \DiscRatio_n^{\text{full}}\}$, batch size $b \leq p$ and aggregation length $L$. 
    \For{$n = 0, 1, \dots$}
        \State Randomly select batch $I_n = \{i_1, \dots, i_b\} \subseteq \{1, \ldots, p\}$ with $|I_n| = b$ via uniform distribution. 
        \State $\rho_{i,n} \gets$ Approximation of $\max\{c_F,\|F_{f_i}'(\sigmad_n)\|_{\mathcal{L}(\Hs s \domain,\, \LLdmnd{\partial \domain})}\}^{-1}$ using the power method
        \If{$\DiscRatio_n < \tau^2$}
            \State Terminate iteration. \Comment{Take $\sigmad_n$ as approximate solution.}
        \EndIf
        \State $\sigmad_{n+1} \gets \sigmad_n - \dfrac{1}{b} \sum_{i \in I_n} \rho_{i,n}^2\, F_{f_i}'(\sigmad_n)^\ast \left(F_{f_i}(\sigma_n) - y_i^\delta\right)$
    \EndFor
\end{algorithmic}
\end{algorithm}
\subsubsection*{Simulation details}
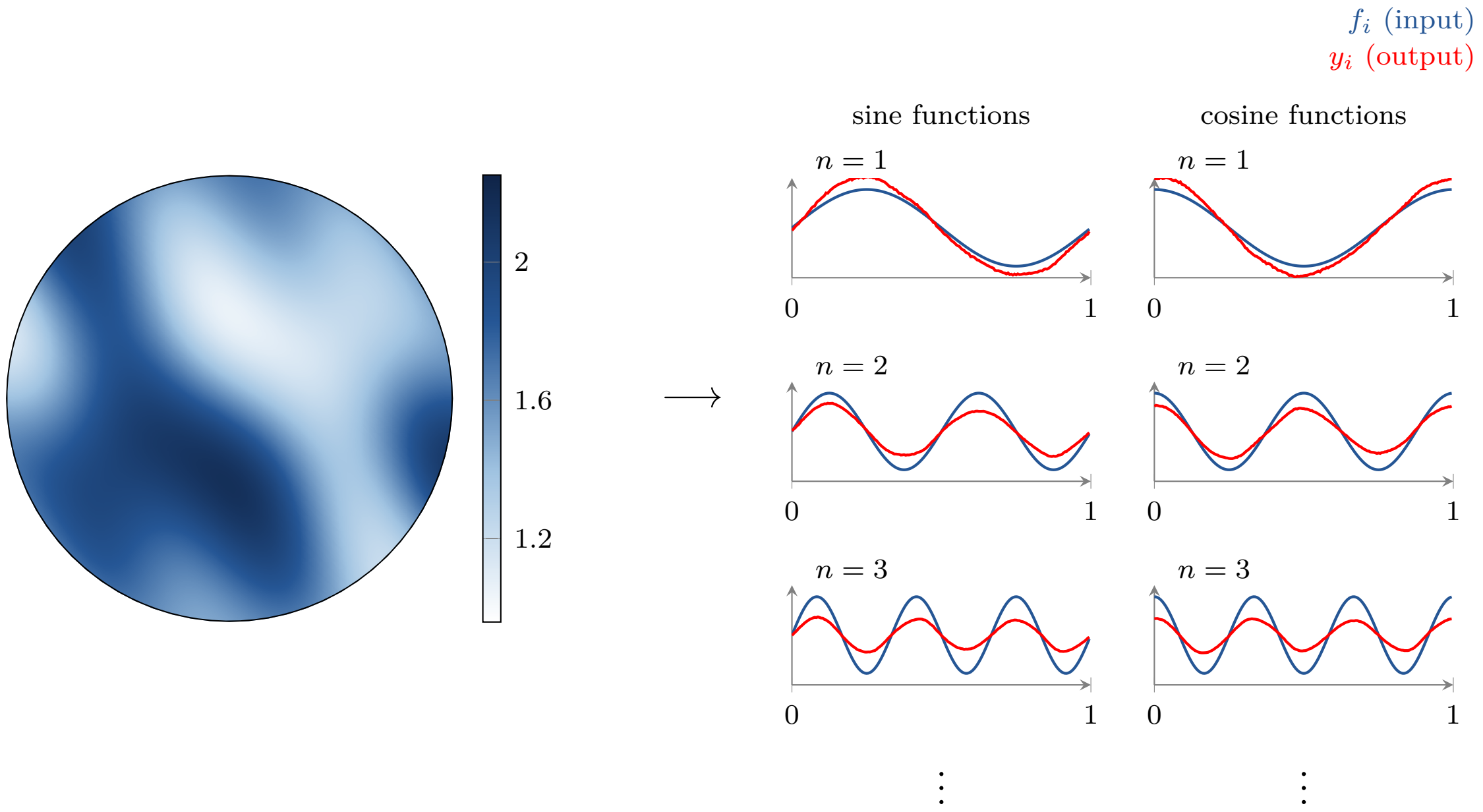
\begin{figure}[ht]
    \centering\tikzset{external/export= true}
    \pgfplotsset{
	compact plot/.style={
		height=1.0cm,width=3cm,scale only axis,
		xmin=0,xmax=1,ymin=-1.3,ymax=1.3,
		axis lines=left,axis line style={draw=gray},
		xtick={0,1},ytick=\empty,tick label style={font=\scriptsize},
		tick style={gray},
		title style={font=\scriptsize,at={(0,1)},anchor=north west,yshift=-0.05cm,xshift=0.1cm},
	}
}
\begin{minipage}[c]{1\textwidth}
\begin{minipage}[t]{0.45\textwidth}
	\centering
    \vspace{3.8em}
	\begin{tikzpicture}
		\begin{axis}[
			width=5.2cm,scale only axis,
			axis equal image,hide axis,
			xmin=-2,xmax=2,ymin=-2,ymax=2,
			colormap name=wuemap,
			view={0}{90},
			colorbar,
			colorbar style={width=0.18cm,ytick={1.2,1.6,2.0},font=\scriptsize},
			]
			\def\circlepath{}
			\foreach \ang in {0,5,...,355} {
				\pgfmathsetmacro{\cx}{2*cos(\ang)}
				\pgfmathsetmacro{\cy}{2*sin(\ang)}
				\xdef\circlepath{\circlepath (axis cs:\cx,\cy) --}
			}
			\path[clip] \circlepath cycle;
			\addplot3[surf,shader=interp,mesh/ordering=y varies] table {figures/tikz/data/sigma_grid.dat};
			\addplot[domain=0:360,samples=100,thick,black] ({2*cos(x)},{2*sin(x)});
		\end{axis}
	\end{tikzpicture}
\end{minipage}\begin{minipage}[t]{0.06\textwidth}
	\centering
    \vspace{9.5em}
	$\longrightarrow$
\end{minipage}\begin{minipage}[t]{0.48\textwidth}
	\centering
	{\small }\hfill{\scriptsize\color{wue}$f_i$ (input)\quad}\\[-0.3em]\hfill{\scriptsize\color{red}$y_i$ (output)\quad}\\[0.4em]
	\begin{tabular}{@{}c@{\hspace{0.6em}}c@{}}
		{\scriptsize sine functions} & {\scriptsize cosine functions}\\
		\begin{tikzpicture}
			\begin{axis}[compact plot,title={$n=1$},
             title style={at={(0,1.02)},anchor=west,xshift=0pt}]
				\addplot[wue,thick] table[x=t,y=sin1] {figures/tikz/data/io_functions.dat};
				\addplot[red,thick] table[x=t,y=usin1] {figures/tikz/data/io_functions.dat};
			\end{axis}
		\end{tikzpicture} &
		\begin{tikzpicture}
			\begin{axis}[compact plot,title={$n=1$}, title style={at={(0,1.02)}, anchor=west}]
				\addplot[wue,thick] table[x=t,y=cos1] {figures/tikz/data/io_functions.dat};
				\addplot[red,thick] table[x=t,y=ucos1] {figures/tikz/data/io_functions.dat};
			\end{axis}
		\end{tikzpicture}\\
		\begin{tikzpicture}
			\begin{axis}[compact plot,title={$n=2$},
            title style={at={(0,1.0)},anchor=west,xshift=0pt}]
				\addplot[wue,thick] table[x=t,y=sin2] {figures/tikz/data/io_functions.dat};
				\addplot[red,thick] table[x=t,y=usin2] {figures/tikz/data/io_functions.dat};
			\end{axis}
		\end{tikzpicture} &
		\begin{tikzpicture}
			\begin{axis}[compact plot,title={$n=2$},
            title style={at={(0,1.0)},anchor=west,xshift=0pt}]
				\addplot[wue,thick] table[x=t,y=cos2] {figures/tikz/data/io_functions.dat};
				\addplot[red,thick] table[x=t,y=ucos2] {figures/tikz/data/io_functions.dat};
			\end{axis}
		\end{tikzpicture}\\
		\begin{tikzpicture}
			\begin{axis}[compact plot,title={$n=3$},
            title style={at={(0,1.0)},anchor=west,xshift=0pt}]
				\addplot[wue,thick] table[x=t,y=sin3] {figures/tikz/data/io_functions.dat};
				\addplot[red,thick] table[x=t,y=usin3] {figures/tikz/data/io_functions.dat};
			\end{axis}
		\end{tikzpicture} &
		\begin{tikzpicture}
			\begin{axis}[compact plot,title={$n=3$},
            title style={at={(0,1.0)},anchor=west,xshift=0pt}]
				\addplot[wue,thick] table[x=t,y=cos3] {figures/tikz/data/io_functions.dat};
				\addplot[red,thick] table[x=t,y=ucos3] {figures/tikz/data/io_functions.dat};
			\end{axis}
		\end{tikzpicture}\\
		$\vdots$ & $\vdots$
	\end{tabular}
\end{minipage}
\end{minipage}     \caption{Example setting: exact conductivity $\bar \sigma$ (left) and selected Neumann input functions $f_i$ with the corresponding boundary data $y_i$ (right). Shown are 6 of the $p=182$ input/output pairs.}
    \label{fig:new_num_example_setting}
\end{figure}

 We consider a two-dimensional test function defined as  
\begin{align*}
    \bar \sigma(x,y) =  1.605&-0.05\sin( 2 x)-0.25\sin( 2 y) -0.25 \sin(1.81x+0.85y)\\+&0.08\sin( 4 x)-0.025\sin( 4 y) 
    -0.15 \sin(3.24x+2.35y)
\end{align*}
on the domain $\domain= \{ (x,y) : x^2 +y^2 \leq 4 \}$ (see the left panel of Figure \ref{fig:new_num_example_setting}). The function is smooth and hence $\bar\sigma \in \Hs s \domain$ for $s>1$. In the experiments we will work with $s:= 1.01$. 
The domain is discretized by an equidistant grid with $64\times 64$ grid points and the finite element method is based on a triangular mesh with 792 degree of freedom.\\
The mesh implies a grid on the boundary with $184$ grid points. Hence as Neumann input functions we select N=182 basis functions of $\LLdmnd{\partial \domain}$, i.e.
\begin{equation*}
 \mathcal B_{\partial \domain}^N :=\left\{ b\in \LLdmnd{\partial \domain}: b(t)= \sin( n t )  \text{ or } b(t)= \cos( n t ) \text{ with } n=1,\ldots,\lfloor N/2 \rfloor  \right\}
\end{equation*}
where $t:\partial \domain \to [0,2\pi)$ corresponds to a uniform parametrization of the boundary (see the right panel of Figure \ref{fig:new_num_example_setting}).\\
The experiments are applied on simulated data. First the discretized forward operator is evaluated for the exact conductivity and then a Gaussian noise is added to the data. In the simulation, the model error is neglected by using the same discretization for both simulation and iteration. This allows the use of high-frequency basis functions, in which the model error dominates.\\
To test the Landweber-Kaczmarz method, we start with the constant function $\sigma_0 \equiv 1$,
get an iteration $\sigmad_n$ and terminate either if the full or batched discrepancy principle holds for the parameter $ \tau= 1.05$. Since all other parameters are fixed, the iteration $\sigmad_n$ depends only on the batch size $b$. However, due to the randomness in the selection of
the batch, the same parameter choice will not produce the same results. Hence, the
experiments are carried out multiple times. Note that due to the Gaussian noise, here with a noise level of $5\%$, the
data also has a random component. This randomness is not addressed in the experiments; instead the
same data is used for each run.
\subsubsection*{Results}
\begin{figure}[ht]
    \centering\tikzset{external/export= false}
    \definecolor{clr0}{rgb}{0.0330,0.1267,0.1916}
\definecolor{clr1}{rgb}{0.0417,0.1600,0.2420}
\definecolor{clr2}{rgb}{0.0504,0.1933,0.2924}
\definecolor{clr3}{rgb}{0.0590,0.2267,0.3429}
\definecolor{clr4}{rgb}{0.0677,0.2600,0.3933}
\definecolor{clr5}{rgb}{0.0764,0.2933,0.4437}
\definecolor{clr6}{rgb}{0.0851,0.3267,0.4941}
\definecolor{clr7}{rgb}{0.0938,0.3600,0.5445}
\definecolor{clr8}{rgb}{0.1025,0.3933,0.5950}
\definecolor{clr9}{rgb}{0.1111,0.4267,0.6454}
\definecolor{clr10}{rgb}{0.1198,0.4600,0.6958}
\definecolor{clr11}{rgb}{0.1285,0.4933,0.7462}
\definecolor{clr12}{rgb}{0.1372,0.5267,0.7966}
\definecolor{clr13}{rgb}{0.1459,0.5600,0.8471}
\definecolor{clr14}{rgb}{0.1546,0.5933,0.8975}
\definecolor{clr15}{rgb}{0.1632,0.6267,0.9479}
\definecolor{clr16}{rgb}{0.1719,0.6600,0.9983}
\definecolor{clr17}{rgb}{0.1806,0.6933,1.0000}
\definecolor{clr18}{rgb}{0.1893,0.7267,1.0000}
\definecolor{clr19}{rgb}{0.1980,0.7600,1.0000}
\definecolor{clr20}{rgb}{0.2067,0.7933,1.0000}
\definecolor{clr21}{rgb}{0.2154,0.8267,1.0000}
\definecolor{clr22}{rgb}{0.2240,0.8600,1.0000}
\definecolor{clr23}{rgb}{0.2327,0.8933,1.0000}
\definecolor{clr24}{rgb}{0.2414,0.9267,1.0000}
\definecolor{clr25}{rgb}{0.2501,0.9600,1.0000}
\definecolor{clr26}{rgb}{0.2588,0.9933,1.0000}
\definecolor{clr27}{rgb}{0.2714,0.1352,0.0149}
\definecolor{clr28}{rgb}{0.3429,0.1708,0.0188}
\definecolor{clr29}{rgb}{0.4143,0.2063,0.0227}
\definecolor{clr30}{rgb}{0.4857,0.2419,0.0267}
\definecolor{clr31}{rgb}{0.5571,0.2775,0.0306}
\definecolor{clr32}{rgb}{0.6286,0.3131,0.0345}
\definecolor{clr33}{rgb}{0.7000,0.3486,0.0384}
\definecolor{clr34}{rgb}{0.7714,0.3842,0.0424}
\definecolor{clr35}{rgb}{0.8429,0.4198,0.0463}
\definecolor{clr36}{rgb}{0.9143,0.4554,0.0502}
\definecolor{clr37}{rgb}{0.9857,0.4909,0.0541}
\definecolor{clr38}{rgb}{1.0000,0.5265,0.0580}
\definecolor{clr39}{rgb}{1.0000,0.5621,0.0620}
\definecolor{clr40}{rgb}{1.0000,0.5976,0.0659}
\definecolor{clr41}{rgb}{1.0000,0.6332,0.0698}
\definecolor{clr42}{rgb}{1.0000,0.6688,0.0737}
\definecolor{clr43}{rgb}{1.0000,0.7044,0.0776}
\definecolor{clr44}{rgb}{1.0000,0.7399,0.0816}
\definecolor{clr45}{rgb}{1.0000,0.7755,0.0855}
\definecolor{clr46}{rgb}{1.0000,0.8111,0.0894}
\definecolor{clr47}{rgb}{1.0000,0.8467,0.0933}
\definecolor{clr48}{rgb}{1.0000,0.8822,0.0973}
\definecolor{clr49}{rgb}{1.0000,0.9178,0.1012}
\definecolor{clr50}{rgb}{1.0000,0.9534,0.1051}
\definecolor{clr51}{rgb}{1.0000,0.9890,0.1090}
\definecolor{clr52}{rgb}{1.0000,1.0000,0.1129}
\definecolor{clr53}{rgb}{1.0000,1.0000,0.1169}
\definecolor{clr54}{rgb}{0.0468,0.1703,0.0468}
\definecolor{clr55}{rgb}{0.0592,0.2151,0.0592}
\definecolor{clr56}{rgb}{0.0715,0.2599,0.0715}
\definecolor{clr57}{rgb}{0.0838,0.3048,0.0838}
\definecolor{clr58}{rgb}{0.0961,0.3496,0.0961}
\definecolor{clr59}{rgb}{0.1085,0.3944,0.1085}
\definecolor{clr60}{rgb}{0.1208,0.4392,0.1208}
\definecolor{clr61}{rgb}{0.1331,0.4840,0.1331}
\definecolor{clr62}{rgb}{0.1454,0.5289,0.1454}
\definecolor{clr63}{rgb}{0.1578,0.5737,0.1578}
\definecolor{clr64}{rgb}{0.1701,0.6185,0.1701}
\definecolor{clr65}{rgb}{0.1824,0.6633,0.1824}
\definecolor{clr66}{rgb}{0.1947,0.7081,0.1947}
\definecolor{clr67}{rgb}{0.2071,0.7529,0.2071}
\definecolor{clr68}{rgb}{0.2194,0.7978,0.2194}
\definecolor{clr69}{rgb}{0.2317,0.8426,0.2317}
\definecolor{clr70}{rgb}{0.2440,0.8874,0.2440}
\definecolor{clr71}{rgb}{0.2564,0.9322,0.2564}
\definecolor{clr72}{rgb}{0.2687,0.9770,0.2687}
\definecolor{clr73}{rgb}{0.2810,1.0000,0.2810}
\definecolor{clr74}{rgb}{0.2933,1.0000,0.2933}
\definecolor{clr75}{rgb}{0.3057,1.0000,0.3057}
\definecolor{clr76}{rgb}{0.3180,1.0000,0.3180}
\definecolor{clr77}{rgb}{0.3303,1.0000,0.3303}
\definecolor{clr78}{rgb}{0.3426,1.0000,0.3426}
\definecolor{clr79}{rgb}{0.3550,1.0000,0.3550}
\definecolor{clr80}{rgb}{0.3673,1.0000,0.3673}
\definecolor{clr81}{rgb}{0.2278,0.0415,0.0426}
\definecolor{clr82}{rgb}{0.2877,0.0524,0.0538}
\definecolor{clr83}{rgb}{0.3477,0.0634,0.0650}
\definecolor{clr84}{rgb}{0.4076,0.0743,0.0762}
\definecolor{clr85}{rgb}{0.4676,0.0852,0.0874}
\definecolor{clr86}{rgb}{0.5275,0.0961,0.0986}
\definecolor{clr87}{rgb}{0.5875,0.1071,0.1098}
\definecolor{clr88}{rgb}{0.6474,0.1180,0.1210}
\definecolor{clr89}{rgb}{0.7073,0.1289,0.1322}
\definecolor{clr90}{rgb}{0.7673,0.1398,0.1434}
\definecolor{clr91}{rgb}{0.8272,0.1508,0.1546}
\definecolor{clr92}{rgb}{0.8872,0.1617,0.1658}
\definecolor{clr93}{rgb}{0.9471,0.1726,0.1770}
\definecolor{clr94}{rgb}{1.0000,0.1835,0.1882}
\definecolor{clr95}{rgb}{1.0000,0.1945,0.1994}
\definecolor{clr96}{rgb}{1.0000,0.2054,0.2106}
\definecolor{clr97}{rgb}{1.0000,0.2163,0.2218}
\definecolor{clr98}{rgb}{1.0000,0.2272,0.2331}
\definecolor{clr99}{rgb}{1.0000,0.2382,0.2443}
\definecolor{clr100}{rgb}{1.0000,0.2491,0.2555}
\definecolor{clr101}{rgb}{1.0000,0.2600,0.2667}
\definecolor{clr102}{rgb}{1.0000,0.2709,0.2779}
\definecolor{clr103}{rgb}{1.0000,0.2818,0.2891}
\definecolor{clr104}{rgb}{1.0000,0.2928,0.3003}
\definecolor{clr105}{rgb}{1.0000,0.3037,0.3115}
\definecolor{clr106}{rgb}{1.0000,0.3146,0.3227}
\definecolor{clr107}{rgb}{1.0000,0.3255,0.3339}
\definecolor{clr108}{rgb}{0.1575,0.1096,0.2012}
\definecolor{clr109}{rgb}{0.1990,0.1385,0.2541}
\definecolor{clr110}{rgb}{0.2404,0.1673,0.3071}
\definecolor{clr111}{rgb}{0.2819,0.1962,0.3600}
\definecolor{clr112}{rgb}{0.3234,0.2250,0.4129}
\definecolor{clr113}{rgb}{0.3648,0.2539,0.4659}
\definecolor{clr114}{rgb}{0.4063,0.2827,0.5188}
\definecolor{clr115}{rgb}{0.4477,0.3116,0.5718}
\definecolor{clr116}{rgb}{0.4892,0.3404,0.6247}
\definecolor{clr117}{rgb}{0.5306,0.3693,0.6776}
\definecolor{clr118}{rgb}{0.5721,0.3982,0.7306}
\definecolor{clr119}{rgb}{0.6136,0.4270,0.7835}
\definecolor{clr120}{rgb}{0.6550,0.4559,0.8365}
\definecolor{clr121}{rgb}{0.6965,0.4847,0.8894}
\definecolor{clr122}{rgb}{0.7379,0.5136,0.9424}
\definecolor{clr123}{rgb}{0.7794,0.5424,0.9953}
\definecolor{clr124}{rgb}{0.8208,0.5713,1.0000}
\definecolor{clr125}{rgb}{0.8623,0.6001,1.0000}
\definecolor{clr126}{rgb}{0.9038,0.6290,1.0000}
\definecolor{clr127}{rgb}{0.9452,0.6578,1.0000}
\definecolor{clr128}{rgb}{0.9867,0.6867,1.0000}
\definecolor{clr129}{rgb}{1.0000,0.7155,1.0000}
\definecolor{clr130}{rgb}{1.0000,0.7444,1.0000}
\definecolor{clr131}{rgb}{1.0000,0.7732,1.0000}
\definecolor{clr132}{rgb}{1.0000,0.8021,1.0000}
\definecolor{clr133}{rgb}{1.0000,0.8309,1.0000}
\definecolor{clr134}{rgb}{1.0000,0.8598,1.0000}
\definecolor{clr135}{rgb}{0.1216,0.4667,0.7059}
\definecolor{clr136}{rgb}{1.0000,0.4980,0.0549}
\definecolor{clr137}{rgb}{0.1725,0.6275,0.1725}
\definecolor{clr138}{rgb}{0.8392,0.1529,0.1569}
\definecolor{clr139}{rgb}{0.5804,0.4039,0.7412}
\definecolor{clr140}{rgb}{0.5490,0.3373,0.2941}
 \pgfplotsset{final_error_yticks/.style={ytick={0.04,0.06,0.08,0.1,0.2},yticklabels={0.04,0.06,0.08,0.1,0.2}}}
 \pgfplotsset{
	final error axis/.style={
		width=9.2cm,height=6.2cm,scale only axis,
		xmin=160,xmax=4e6,ymin=0.036,ymax=0.2,
		at={(0,0)},anchor=origin,
	},
}
\begin{tikzpicture}
	\begin{axis}[
		final error axis,xmode=log,ymode=log,
		name=mainax,
		xlabel={Number of PDE-solves},
		ylabel={Relative $L^2$-error},
		label style={font=\small},
		tick label style={font=\footnotesize},
		final_error_yticks,
		grid=major,
		grid style={gray!25},
		axis background/.style={fill=gray!8},
		]
		\addplot[draw=none] coordinates {(1000,0.1)};
	\end{axis}
	\begin{axis}[final error axis,xmode=log,ymode=log,hide axis]
		\addplot[only marks, mark=square*, mark size=0.45pt, color=clr0, mark options={draw=none}] table {figures/tikz/data/final_error_b1_L1.dat};
		\addplot[only marks, mark=square*, mark size=0.45pt, color=clr1, mark options={draw=none}] table {figures/tikz/data/final_error_b1_L2.dat};
		\addplot[only marks, mark=square*, mark size=0.45pt, color=clr2, mark options={draw=none}] table {figures/tikz/data/final_error_b1_L3.dat};
		\addplot[only marks, mark=square*, mark size=0.45pt, color=clr3, mark options={draw=none}] table {figures/tikz/data/final_error_b1_L4.dat};
		\addplot[only marks, mark=square*, mark size=0.45pt, color=clr4, mark options={draw=none}] table {figures/tikz/data/final_error_b1_L5.dat};
		\addplot[only marks, mark=square*, mark size=0.45pt, color=clr5, mark options={draw=none}] table {figures/tikz/data/final_error_b1_L6.dat};
		\addplot[only marks, mark=square*, mark size=0.45pt, color=clr6, mark options={draw=none}] table {figures/tikz/data/final_error_b1_L7.dat};
		\addplot[only marks, mark=square*, mark size=0.45pt, color=clr7, mark options={draw=none}] table {figures/tikz/data/final_error_b1_L8.dat};
		\addplot[only marks, mark=square*, mark size=0.45pt, color=clr8, mark options={draw=none}] table {figures/tikz/data/final_error_b1_L9.dat};
		\addplot[only marks, mark=square*, mark size=0.45pt, color=clr9, mark options={draw=none}] table {figures/tikz/data/final_error_b1_L10.dat};
		\addplot[only marks, mark=square*, mark size=0.45pt, color=clr10, mark options={draw=none}] table {figures/tikz/data/final_error_b1_L20.dat};
		\addplot[only marks, mark=square*, mark size=0.45pt, color=clr11, mark options={draw=none}] table {figures/tikz/data/final_error_b1_L30.dat};
		\addplot[only marks, mark=square*, mark size=0.45pt, color=clr12, mark options={draw=none}] table {figures/tikz/data/final_error_b1_L40.dat};
		\addplot[only marks, mark=square*, mark size=0.45pt, color=clr13, mark options={draw=none}] table {figures/tikz/data/final_error_b1_L50.dat};
		\addplot[only marks, mark=square*, mark size=0.45pt, color=clr14, mark options={draw=none}] table {figures/tikz/data/final_error_b1_L60.dat};
		\addplot[only marks, mark=square*, mark size=0.45pt, color=clr15, mark options={draw=none}] table {figures/tikz/data/final_error_b1_L70.dat};
	\end{axis}
	\begin{axis}[final error axis,xmode=log,ymode=log,hide axis]
		\addplot[only marks, mark=square*, mark size=0.45pt, color=clr16, mark options={draw=none}] table {figures/tikz/data/final_error_b1_L80.dat};
		\addplot[only marks, mark=square*, mark size=0.45pt, color=clr17, mark options={draw=none}] table {figures/tikz/data/final_error_b1_L90.dat};
		\addplot[only marks, mark=square*, mark size=0.45pt, color=clr18, mark options={draw=none}] table {figures/tikz/data/final_error_b1_L100.dat};
		\addplot[only marks, mark=square*, mark size=0.45pt, color=clr19, mark options={draw=none}] table {figures/tikz/data/final_error_b1_L150.dat};
		\addplot[only marks, mark=square*, mark size=0.45pt, color=clr20, mark options={draw=none}] table {figures/tikz/data/final_error_b1_L200.dat};
		\addplot[only marks, mark=square*, mark size=0.45pt, color=clr21, mark options={draw=none}] table {figures/tikz/data/final_error_b1_L250.dat};
		\addplot[only marks, mark=square*, mark size=0.45pt, color=clr22, mark options={draw=none}] table {figures/tikz/data/final_error_b1_L300.dat};
		\addplot[only marks, mark=square*, mark size=0.45pt, color=clr23, mark options={draw=none}] table {figures/tikz/data/final_error_b1_L350.dat};
		\addplot[only marks, mark=square*, mark size=0.45pt, color=clr24, mark options={draw=none}] table {figures/tikz/data/final_error_b1_L400.dat};
		\addplot[only marks, mark=square*, mark size=0.45pt, color=clr25, mark options={draw=none}] table {figures/tikz/data/final_error_b1_L450.dat};
		\addplot[only marks, mark=square*, mark size=0.45pt, color=clr26, mark options={draw=none}] table {figures/tikz/data/final_error_b1_L500.dat};
		\addplot[only marks, mark=square*, mark size=0.45pt, color=clr27, mark options={draw=none}] table {figures/tikz/data/final_error_b3_L1.dat};
		\addplot[only marks, mark=square*, mark size=0.45pt, color=clr28, mark options={draw=none}] table {figures/tikz/data/final_error_b3_L2.dat};
		\addplot[only marks, mark=square*, mark size=0.45pt, color=clr29, mark options={draw=none}] table {figures/tikz/data/final_error_b3_L3.dat};
		\addplot[only marks, mark=square*, mark size=0.45pt, color=clr30, mark options={draw=none}] table {figures/tikz/data/final_error_b3_L4.dat};
		\addplot[only marks, mark=square*, mark size=0.45pt, color=clr31, mark options={draw=none}] table {figures/tikz/data/final_error_b3_L5.dat};
	\end{axis}
	\begin{axis}[final error axis,xmode=log,ymode=log,hide axis]
		\addplot[only marks, mark=square*, mark size=0.45pt, color=clr32, mark options={draw=none}] table {figures/tikz/data/final_error_b3_L6.dat};
		\addplot[only marks, mark=square*, mark size=0.45pt, color=clr33, mark options={draw=none}] table {figures/tikz/data/final_error_b3_L7.dat};
		\addplot[only marks, mark=square*, mark size=0.45pt, color=clr34, mark options={draw=none}] table {figures/tikz/data/final_error_b3_L8.dat};
		\addplot[only marks, mark=square*, mark size=0.45pt, color=clr35, mark options={draw=none}] table {figures/tikz/data/final_error_b3_L9.dat};
		\addplot[only marks, mark=square*, mark size=0.45pt, color=clr36, mark options={draw=none}] table {figures/tikz/data/final_error_b3_L10.dat};
		\addplot[only marks, mark=square*, mark size=0.45pt, color=clr37, mark options={draw=none}] table {figures/tikz/data/final_error_b3_L20.dat};
		\addplot[only marks, mark=square*, mark size=0.45pt, color=clr38, mark options={draw=none}] table {figures/tikz/data/final_error_b3_L30.dat};
		\addplot[only marks, mark=square*, mark size=0.45pt, color=clr39, mark options={draw=none}] table {figures/tikz/data/final_error_b3_L40.dat};
		\addplot[only marks, mark=square*, mark size=0.45pt, color=clr40, mark options={draw=none}] table {figures/tikz/data/final_error_b3_L50.dat};
		\addplot[only marks, mark=square*, mark size=0.45pt, color=clr41, mark options={draw=none}] table {figures/tikz/data/final_error_b3_L60.dat};
		\addplot[only marks, mark=square*, mark size=0.45pt, color=clr42, mark options={draw=none}] table {figures/tikz/data/final_error_b3_L70.dat};
		\addplot[only marks, mark=square*, mark size=0.45pt, color=clr43, mark options={draw=none}] table {figures/tikz/data/final_error_b3_L80.dat};
		\addplot[only marks, mark=square*, mark size=0.45pt, color=clr44, mark options={draw=none}] table {figures/tikz/data/final_error_b3_L90.dat};
		\addplot[only marks, mark=square*, mark size=0.45pt, color=clr45, mark options={draw=none}] table {figures/tikz/data/final_error_b3_L100.dat};
		\addplot[only marks, mark=square*, mark size=0.45pt, color=clr46, mark options={draw=none}] table {figures/tikz/data/final_error_b3_L150.dat};
		\addplot[only marks, mark=square*, mark size=0.45pt, color=clr47, mark options={draw=none}] table {figures/tikz/data/final_error_b3_L200.dat};
	\end{axis}
	\begin{axis}[final error axis,xmode=log,ymode=log,hide axis]
		\addplot[only marks, mark=square*, mark size=0.45pt, color=clr48, mark options={draw=none}] table {figures/tikz/data/final_error_b3_L250.dat};
		\addplot[only marks, mark=square*, mark size=0.45pt, color=clr49, mark options={draw=none}] table {figures/tikz/data/final_error_b3_L300.dat};
		\addplot[only marks, mark=square*, mark size=0.45pt, color=clr50, mark options={draw=none}] table {figures/tikz/data/final_error_b3_L350.dat};
		\addplot[only marks, mark=square*, mark size=0.45pt, color=clr51, mark options={draw=none}] table {figures/tikz/data/final_error_b3_L400.dat};
		\addplot[only marks, mark=square*, mark size=0.45pt, color=clr52, mark options={draw=none}] table {figures/tikz/data/final_error_b3_L450.dat};
		\addplot[only marks, mark=square*, mark size=0.45pt, color=clr53, mark options={draw=none}] table {figures/tikz/data/final_error_b3_L500.dat};
		\addplot[only marks, mark=square*, mark size=0.45pt, color=clr54, mark options={draw=none}] table {figures/tikz/data/final_error_b10_L1.dat};
		\addplot[only marks, mark=square*, mark size=0.45pt, color=clr55, mark options={draw=none}] table {figures/tikz/data/final_error_b10_L2.dat};
		\addplot[only marks, mark=square*, mark size=0.45pt, color=clr56, mark options={draw=none}] table {figures/tikz/data/final_error_b10_L3.dat};
		\addplot[only marks, mark=square*, mark size=0.45pt, color=clr57, mark options={draw=none}] table {figures/tikz/data/final_error_b10_L4.dat};
		\addplot[only marks, mark=square*, mark size=0.45pt, color=clr58, mark options={draw=none}] table {figures/tikz/data/final_error_b10_L5.dat};
		\addplot[only marks, mark=square*, mark size=0.45pt, color=clr59, mark options={draw=none}] table {figures/tikz/data/final_error_b10_L6.dat};
		\addplot[only marks, mark=square*, mark size=0.45pt, color=clr60, mark options={draw=none}] table {figures/tikz/data/final_error_b10_L7.dat};
		\addplot[only marks, mark=square*, mark size=0.45pt, color=clr61, mark options={draw=none}] table {figures/tikz/data/final_error_b10_L8.dat};
		\addplot[only marks, mark=square*, mark size=0.45pt, color=clr62, mark options={draw=none}] table {figures/tikz/data/final_error_b10_L9.dat};
		\addplot[only marks, mark=square*, mark size=0.45pt, color=clr63, mark options={draw=none}] table {figures/tikz/data/final_error_b10_L10.dat};
	\end{axis}
	\begin{axis}[final error axis,xmode=log,ymode=log,hide axis]
		\addplot[only marks, mark=square*, mark size=0.45pt, color=clr64, mark options={draw=none}] table {figures/tikz/data/final_error_b10_L20.dat};
		\addplot[only marks, mark=square*, mark size=0.45pt, color=clr65, mark options={draw=none}] table {figures/tikz/data/final_error_b10_L30.dat};
		\addplot[only marks, mark=square*, mark size=0.45pt, color=clr66, mark options={draw=none}] table {figures/tikz/data/final_error_b10_L40.dat};
		\addplot[only marks, mark=square*, mark size=0.45pt, color=clr67, mark options={draw=none}] table {figures/tikz/data/final_error_b10_L50.dat};
		\addplot[only marks, mark=square*, mark size=0.45pt, color=clr68, mark options={draw=none}] table {figures/tikz/data/final_error_b10_L60.dat};
		\addplot[only marks, mark=square*, mark size=0.45pt, color=clr69, mark options={draw=none}] table {figures/tikz/data/final_error_b10_L70.dat};
		\addplot[only marks, mark=square*, mark size=0.45pt, color=clr70, mark options={draw=none}] table {figures/tikz/data/final_error_b10_L80.dat};
		\addplot[only marks, mark=square*, mark size=0.45pt, color=clr71, mark options={draw=none}] table {figures/tikz/data/final_error_b10_L90.dat};
		\addplot[only marks, mark=square*, mark size=0.45pt, color=clr72, mark options={draw=none}] table {figures/tikz/data/final_error_b10_L100.dat};
		\addplot[only marks, mark=square*, mark size=0.45pt, color=clr73, mark options={draw=none}] table {figures/tikz/data/final_error_b10_L150.dat};
		\addplot[only marks, mark=square*, mark size=0.45pt, color=clr74, mark options={draw=none}] table {figures/tikz/data/final_error_b10_L200.dat};
		\addplot[only marks, mark=square*, mark size=0.45pt, color=clr75, mark options={draw=none}] table {figures/tikz/data/final_error_b10_L250.dat};
		\addplot[only marks, mark=square*, mark size=0.45pt, color=clr76, mark options={draw=none}] table {figures/tikz/data/final_error_b10_L300.dat};
		\addplot[only marks, mark=square*, mark size=0.45pt, color=clr77, mark options={draw=none}] table {figures/tikz/data/final_error_b10_L350.dat};
		\addplot[only marks, mark=square*, mark size=0.45pt, color=clr78, mark options={draw=none}] table {figures/tikz/data/final_error_b10_L400.dat};
		\addplot[only marks, mark=square*, mark size=0.45pt, color=clr79, mark options={draw=none}] table {figures/tikz/data/final_error_b10_L450.dat};
	\end{axis}
	\begin{axis}[final error axis,xmode=log,ymode=log,hide axis]
		\addplot[only marks, mark=square*, mark size=0.45pt, color=clr80, mark options={draw=none}] table {figures/tikz/data/final_error_b10_L500.dat};
		\addplot[only marks, mark=square*, mark size=0.45pt, color=clr81, mark options={draw=none}] table {figures/tikz/data/final_error_b20_L1.dat};
		\addplot[only marks, mark=square*, mark size=0.45pt, color=clr82, mark options={draw=none}] table {figures/tikz/data/final_error_b20_L2.dat};
		\addplot[only marks, mark=square*, mark size=0.45pt, color=clr83, mark options={draw=none}] table {figures/tikz/data/final_error_b20_L3.dat};
		\addplot[only marks, mark=square*, mark size=0.45pt, color=clr84, mark options={draw=none}] table {figures/tikz/data/final_error_b20_L4.dat};
		\addplot[only marks, mark=square*, mark size=0.45pt, color=clr85, mark options={draw=none}] table {figures/tikz/data/final_error_b20_L5.dat};
		\addplot[only marks, mark=square*, mark size=0.45pt, color=clr86, mark options={draw=none}] table {figures/tikz/data/final_error_b20_L6.dat};
		\addplot[only marks, mark=square*, mark size=0.45pt, color=clr87, mark options={draw=none}] table {figures/tikz/data/final_error_b20_L7.dat};
		\addplot[only marks, mark=square*, mark size=0.45pt, color=clr88, mark options={draw=none}] table {figures/tikz/data/final_error_b20_L8.dat};
		\addplot[only marks, mark=square*, mark size=0.45pt, color=clr89, mark options={draw=none}] table {figures/tikz/data/final_error_b20_L9.dat};
		\addplot[only marks, mark=square*, mark size=0.45pt, color=clr90, mark options={draw=none}] table {figures/tikz/data/final_error_b20_L10.dat};
		\addplot[only marks, mark=square*, mark size=0.45pt, color=clr91, mark options={draw=none}] table {figures/tikz/data/final_error_b20_L20.dat};
		\addplot[only marks, mark=square*, mark size=0.45pt, color=clr92, mark options={draw=none}] table {figures/tikz/data/final_error_b20_L30.dat};
		\addplot[only marks, mark=square*, mark size=0.45pt, color=clr93, mark options={draw=none}] table {figures/tikz/data/final_error_b20_L40.dat};
		\addplot[only marks, mark=square*, mark size=0.45pt, color=clr94, mark options={draw=none}] table {figures/tikz/data/final_error_b20_L50.dat};
		\addplot[only marks, mark=square*, mark size=0.45pt, color=clr95, mark options={draw=none}] table {figures/tikz/data/final_error_b20_L60.dat};
	\end{axis}
	\begin{axis}[final error axis,xmode=log,ymode=log,hide axis]
		\addplot[only marks, mark=square*, mark size=0.45pt, color=clr96, mark options={draw=none}] table {figures/tikz/data/final_error_b20_L70.dat};
		\addplot[only marks, mark=square*, mark size=0.45pt, color=clr97, mark options={draw=none}] table {figures/tikz/data/final_error_b20_L80.dat};
		\addplot[only marks, mark=square*, mark size=0.45pt, color=clr98, mark options={draw=none}] table {figures/tikz/data/final_error_b20_L90.dat};
		\addplot[only marks, mark=square*, mark size=0.45pt, color=clr99, mark options={draw=none}] table {figures/tikz/data/final_error_b20_L100.dat};
		\addplot[only marks, mark=square*, mark size=0.45pt, color=clr100, mark options={draw=none}] table {figures/tikz/data/final_error_b20_L150.dat};
		\addplot[only marks, mark=square*, mark size=0.45pt, color=clr101, mark options={draw=none}] table {figures/tikz/data/final_error_b20_L200.dat};
		\addplot[only marks, mark=square*, mark size=0.45pt, color=clr102, mark options={draw=none}] table {figures/tikz/data/final_error_b20_L250.dat};
		\addplot[only marks, mark=square*, mark size=0.45pt, color=clr103, mark options={draw=none}] table {figures/tikz/data/final_error_b20_L300.dat};
		\addplot[only marks, mark=square*, mark size=0.45pt, color=clr104, mark options={draw=none}] table {figures/tikz/data/final_error_b20_L350.dat};
		\addplot[only marks, mark=square*, mark size=0.45pt, color=clr105, mark options={draw=none}] table {figures/tikz/data/final_error_b20_L400.dat};
		\addplot[only marks, mark=square*, mark size=0.45pt, color=clr106, mark options={draw=none}] table {figures/tikz/data/final_error_b20_L450.dat};
		\addplot[only marks, mark=square*, mark size=0.45pt, color=clr107, mark options={draw=none}] table {figures/tikz/data/final_error_b20_L500.dat};
		\addplot[only marks, mark=square*, mark size=0.45pt, color=clr108, mark options={draw=none}] table {figures/tikz/data/final_error_b40_L1.dat};
		\addplot[only marks, mark=square*, mark size=0.45pt, color=clr109, mark options={draw=none}] table {figures/tikz/data/final_error_b40_L2.dat};
		\addplot[only marks, mark=square*, mark size=0.45pt, color=clr110, mark options={draw=none}] table {figures/tikz/data/final_error_b40_L3.dat};
		\addplot[only marks, mark=square*, mark size=0.45pt, color=clr111, mark options={draw=none}] table {figures/tikz/data/final_error_b40_L4.dat};
	\end{axis}
	\begin{axis}[final error axis,xmode=log,ymode=log,hide axis]
		\addplot[only marks, mark=square*, mark size=0.45pt, color=clr112, mark options={draw=none}] table {figures/tikz/data/final_error_b40_L5.dat};
		\addplot[only marks, mark=square*, mark size=0.45pt, color=clr113, mark options={draw=none}] table {figures/tikz/data/final_error_b40_L6.dat};
		\addplot[only marks, mark=square*, mark size=0.45pt, color=clr114, mark options={draw=none}] table {figures/tikz/data/final_error_b40_L7.dat};
		\addplot[only marks, mark=square*, mark size=0.45pt, color=clr115, mark options={draw=none}] table {figures/tikz/data/final_error_b40_L8.dat};
		\addplot[only marks, mark=square*, mark size=0.45pt, color=clr116, mark options={draw=none}] table {figures/tikz/data/final_error_b40_L9.dat};
		\addplot[only marks, mark=square*, mark size=0.45pt, color=clr117, mark options={draw=none}] table {figures/tikz/data/final_error_b40_L10.dat};
		\addplot[only marks, mark=square*, mark size=0.45pt, color=clr118, mark options={draw=none}] table {figures/tikz/data/final_error_b40_L20.dat};
		\addplot[only marks, mark=square*, mark size=0.45pt, color=clr119, mark options={draw=none}] table {figures/tikz/data/final_error_b40_L30.dat};
		\addplot[only marks, mark=square*, mark size=0.45pt, color=clr120, mark options={draw=none}] table {figures/tikz/data/final_error_b40_L40.dat};
		\addplot[only marks, mark=square*, mark size=0.45pt, color=clr121, mark options={draw=none}] table {figures/tikz/data/final_error_b40_L50.dat};
		\addplot[only marks, mark=square*, mark size=0.45pt, color=clr122, mark options={draw=none}] table {figures/tikz/data/final_error_b40_L60.dat};
		\addplot[only marks, mark=square*, mark size=0.45pt, color=clr123, mark options={draw=none}] table {figures/tikz/data/final_error_b40_L70.dat};
		\addplot[only marks, mark=square*, mark size=0.45pt, color=clr124, mark options={draw=none}] table {figures/tikz/data/final_error_b40_L80.dat};
		\addplot[only marks, mark=square*, mark size=0.45pt, color=clr125, mark options={draw=none}] table {figures/tikz/data/final_error_b40_L90.dat};
		\addplot[only marks, mark=square*, mark size=0.45pt, color=clr126, mark options={draw=none}] table {figures/tikz/data/final_error_b40_L100.dat};
		\addplot[only marks, mark=square*, mark size=0.45pt, color=clr127, mark options={draw=none}] table {figures/tikz/data/final_error_b40_L150.dat};
	\end{axis}
	\begin{axis}[final error axis,xmode=log,ymode=log,hide axis]
		\addplot[only marks, mark=square*, mark size=0.45pt, color=clr128, mark options={draw=none}] table {figures/tikz/data/final_error_b40_L200.dat};
		\addplot[only marks, mark=square*, mark size=0.45pt, color=clr129, mark options={draw=none}] table {figures/tikz/data/final_error_b40_L250.dat};
		\addplot[only marks, mark=square*, mark size=0.45pt, color=clr130, mark options={draw=none}] table {figures/tikz/data/final_error_b40_L300.dat};
		\addplot[only marks, mark=square*, mark size=0.45pt, color=clr131, mark options={draw=none}] table {figures/tikz/data/final_error_b40_L350.dat};
		\addplot[only marks, mark=square*, mark size=0.45pt, color=clr132, mark options={draw=none}] table {figures/tikz/data/final_error_b40_L400.dat};
		\addplot[only marks, mark=square*, mark size=0.45pt, color=clr133, mark options={draw=none}] table {figures/tikz/data/final_error_b40_L450.dat};
		\addplot[only marks, mark=square*, mark size=0.45pt, color=clr134, mark options={draw=none}] table {figures/tikz/data/final_error_b40_L500.dat};
		\addplot[only marks, mark=x, mark size=2.00pt, mark options={line width=0.5pt}, draw=clr135] table {figures/tikz/data/final_error_b1_full.dat};
		\addplot[only marks, mark=x, mark size=2.00pt, mark options={line width=0.5pt}, draw=clr136] table {figures/tikz/data/final_error_b3_full.dat};
		\addplot[only marks, mark=x, mark size=2.00pt, mark options={line width=0.5pt}, draw=clr137] table {figures/tikz/data/final_error_b10_full.dat};
		\addplot[only marks, mark=x, mark size=2.00pt, mark options={line width=0.5pt}, draw=clr138] table {figures/tikz/data/final_error_b20_full.dat};
		\addplot[only marks, mark=x, mark size=2.00pt, mark options={line width=0.5pt}, draw=clr139] table {figures/tikz/data/final_error_b40_full.dat};
		\addplot[only marks, mark=x, mark size=2.00pt, mark options={line width=0.5pt}, draw=clr140] table {figures/tikz/data/final_error_b182_full.dat};
	\end{axis}
\node[anchor=north west, inner sep=0pt]
		at ([xshift=4mm]mainax.north east) {\scriptsize
		\vtop{\hbox{\textbf{Discrepancy type}}\vskip 3pt
			\hbox{\begin{tikzpicture}
\draw[black,line width=1.2pt] (0.00,0.00) -- (0.30,0.30) (0.00,0.30) -- (0.30,0.00);
\node[right, font=\tiny] at (0.38,0.15) {classical};
\fill[black] (0.15,-0.25) circle (0.15);
\node[right, font=\tiny] at (0.38,-0.25) {batched};
\end{tikzpicture}
 \unskip}\vskip 10pt
			\hbox{\textbf{Batch size $b$}}\vskip 4pt
			\hbox{\begin{tikzpicture}
\fill[clr135] (0.00,-0.00) rectangle (0.25,0.25);
\node[right, font=\tiny] at (0.32,0.12) {1};
\fill[clr136] (0.00,-0.35) rectangle (0.25,-0.10);
\node[right, font=\tiny] at (0.32,-0.23) {3};
\fill[clr137] (0.00,-0.70) rectangle (0.25,-0.45);
\node[right, font=\tiny] at (0.32,-0.58) {10};
\fill[clr138] (0.95,-0.00) rectangle (1.20,0.25);
\node[right, font=\tiny] at (1.27,0.12) {20};
\fill[clr139] (0.95,-0.35) rectangle (1.20,-0.10);
\node[right, font=\tiny] at (1.27,-0.23) {40};
\fill[clr140] (0.95,-0.70) rectangle (1.20,-0.45);
\node[right, font=\tiny] at (1.27,-0.58) {182};
\end{tikzpicture}
 \unskip}\vskip 10pt
			\hbox{\textbf{Parameter $L$}}\vskip 4pt
			\hbox{\tikzset{every picture/.append style={yscale=0.6}}\begin{tikzpicture}
\fill[black!100.0] (0,0.0000) rectangle (0.25,0.1630);
\fill[black!97.1] (0,0.1630) rectangle (0.25,0.3259);
\fill[black!94.2] (0,0.3259) rectangle (0.25,0.4889);
\fill[black!91.3] (0,0.4889) rectangle (0.25,0.6519);
\fill[black!88.5] (0,0.6519) rectangle (0.25,0.8148);
\fill[black!85.6] (0,0.8148) rectangle (0.25,0.9778);
\fill[black!82.7] (0,0.9778) rectangle (0.25,1.1407);
\fill[black!79.8] (0,1.1407) rectangle (0.25,1.3037);
\fill[black!76.9] (0,1.3037) rectangle (0.25,1.4667);
\fill[black!74.0] (0,1.4667) rectangle (0.25,1.6296);
\fill[black!71.2] (0,1.6296) rectangle (0.25,1.7926);
\fill[black!68.3] (0,1.7926) rectangle (0.25,1.9556);
\fill[black!65.4] (0,1.9556) rectangle (0.25,2.1185);
\fill[black!62.5] (0,2.1185) rectangle (0.25,2.2815);
\fill[black!59.6] (0,2.2815) rectangle (0.25,2.4444);
\fill[black!56.7] (0,2.4444) rectangle (0.25,2.6074);
\fill[black!53.8] (0,2.6074) rectangle (0.25,2.7704);
\fill[black!51.0] (0,2.7704) rectangle (0.25,2.9333);
\fill[black!48.1] (0,2.9333) rectangle (0.25,3.0963);
\fill[black!45.2] (0,3.0963) rectangle (0.25,3.2593);
\fill[black!42.3] (0,3.2593) rectangle (0.25,3.4222);
\fill[black!39.4] (0,3.4222) rectangle (0.25,3.5852);
\fill[black!36.5] (0,3.5852) rectangle (0.25,3.7481);
\fill[black!33.7] (0,3.7481) rectangle (0.25,3.9111);
\fill[black!30.8] (0,3.9111) rectangle (0.25,4.0741);
\fill[black!27.9] (0,4.0741) rectangle (0.25,4.2370);
\fill[black!25.0] (0,4.2370) rectangle (0.25,4.4000);
\draw[line width=0.3pt] (0,0) rectangle (0.25,4.40);
\draw[line width=0.3pt] (0.25,0.0815) -- (0.31,0.0815);
\node[right, font=\tiny, inner sep=1pt] at (0.31,0.0815) {1};
\draw[line width=0.3pt] (0.25,0.7333) -- (0.31,0.7333);
\node[right, font=\tiny, inner sep=1pt] at (0.31,0.7333) {5};
\draw[line width=0.3pt] (0.25,1.5481) -- (0.31,1.5481);
\node[right, font=\tiny, inner sep=1pt] at (0.31,1.5481) {10};
\draw[line width=0.3pt] (0.25,2.2000) -- (0.31,2.2000);
\node[right, font=\tiny, inner sep=1pt] at (0.31,2.2000) {50};
\draw[line width=0.3pt] (0.25,3.0148) -- (0.31,3.0148);
\node[right, font=\tiny, inner sep=1pt] at (0.31,3.0148) {100};
\draw[line width=0.3pt] (0.25,4.3185) -- (0.31,4.3185);
\node[right, font=\tiny, inner sep=1pt] at (0.31,4.3185) {500};
\end{tikzpicture}
 \unskip}
}};
\end{tikzpicture}
     \caption{Relative $L^2$-error of the final iterates versus the number of PDE-solves for the classical and the batched discrepancy principle with different batch sizes and aggregation lengths $L$.
}
    \label{fig:new_num_error_results}
\end{figure}
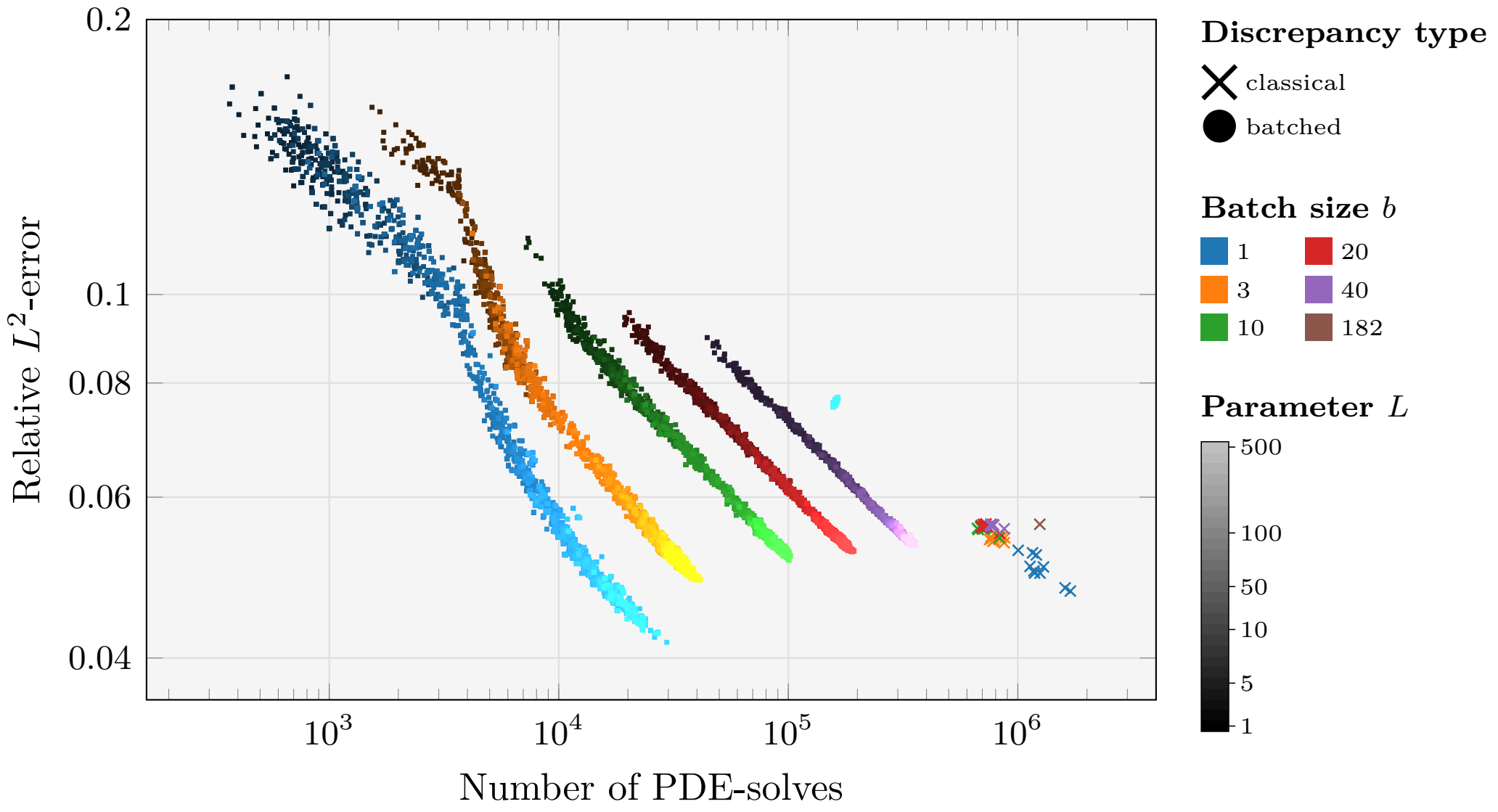

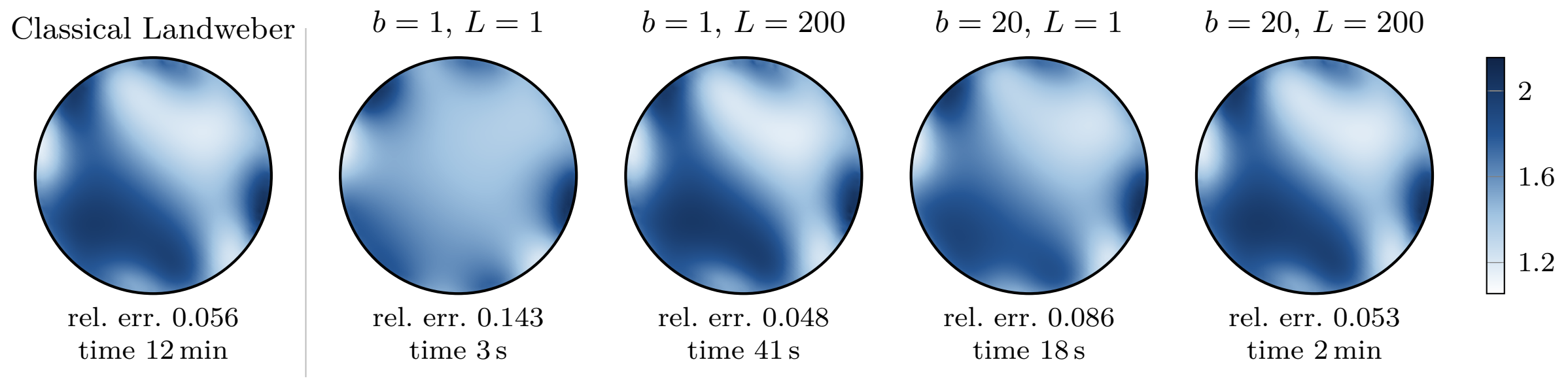
\begin{figure}[ht]
    \centering\tikzset{external/export= false}
    \resizebox{\linewidth}{!}{\begin{tikzpicture}
\begin{axis}[name=recclassical,
  at={(0.00cm,0.00cm)},
  anchor=north west,
  width=2.4cm,
  height=2.4cm,
  scale only axis,
  axis equal image,
  hide axis,
  xmin=-2,xmax=2,ymin=-2,ymax=2,
  view={0}{90},
  colormap name=wuemap,
  point meta min=1.055,
  point meta max=2.155,
  title={Classical Landweber},
  title style={font=\footnotesize,yshift=-0.15cm}]
\clip (axis cs:0,0) circle [radius=2];
\addplot3[surf,shader=interp,mesh/ordering=x varies,mesh/cols=64] table {figures/tikz/data/rec_classical.dat};
\end{axis}
\draw[black,thick] (recclassical.center) circle [radius=1.200cm];
\node[below,font=\scriptsize,align=center] at (recclassical.south) {rel.\ err.\ 0.056\\ time 12\,min};
\begin{axis}[name=recb1L1,
  at={(3.10cm,0.00cm)},
  anchor=north west,
  width=2.4cm,
  height=2.4cm,
  scale only axis,
  axis equal image,
  hide axis,
  xmin=-2,xmax=2,ymin=-2,ymax=2,
  view={0}{90},
  colormap name=wuemap,
  point meta min=1.055,
  point meta max=2.155,
  title={$b=1$, $L=1$},
  title style={font=\footnotesize,yshift=-0.15cm}]
\clip (axis cs:0,0) circle [radius=2];
\addplot3[surf,shader=interp,mesh/ordering=x varies,mesh/cols=64] table {figures/tikz/data/rec_b1_L1.dat};
\end{axis}
\draw[black,thick] (recb1L1.center) circle [radius=1.200cm];
\node[below,font=\scriptsize,align=center] at (recb1L1.south) {rel.\ err.\ 0.143\\ time 3\,s};
\begin{axis}[name=recb1L200,
  at={(6.00cm,0.00cm)},
  anchor=north west,
  width=2.4cm,
  height=2.4cm,
  scale only axis,
  axis equal image,
  hide axis,
  xmin=-2,xmax=2,ymin=-2,ymax=2,
  view={0}{90},
  colormap name=wuemap,
  point meta min=1.055,
  point meta max=2.155,
  title={$b=1$, $L=200$},
  title style={font=\footnotesize,yshift=-0.15cm}]
\clip (axis cs:0,0) circle [radius=2];
\addplot3[surf,shader=interp,mesh/ordering=x varies,mesh/cols=64] table {figures/tikz/data/rec_b1_L200.dat};
\end{axis}
\draw[black,thick] (recb1L200.center) circle [radius=1.200cm];
\node[below,font=\scriptsize,align=center] at (recb1L200.south) {rel.\ err.\ 0.048\\ time 41\,s};
\begin{axis}[name=recb20L1,
  at={(8.90cm,0.00cm)},
  anchor=north west,
  width=2.4cm,
  height=2.4cm,
  scale only axis,
  axis equal image,
  hide axis,
  xmin=-2,xmax=2,ymin=-2,ymax=2,
  view={0}{90},
  colormap name=wuemap,
  point meta min=1.055,
  point meta max=2.155,
  title={$b=20$, $L=1$},
  title style={font=\footnotesize,yshift=-0.15cm}]
\clip (axis cs:0,0) circle [radius=2];
\addplot3[surf,shader=interp,mesh/ordering=x varies,mesh/cols=64] table {figures/tikz/data/rec_b20_L1.dat};
\end{axis}
\draw[black,thick] (recb20L1.center) circle [radius=1.200cm];
\node[below,font=\scriptsize,align=center] at (recb20L1.south) {rel.\ err.\ 0.086\\ time 18\,s};
\begin{axis}[name=recb20L200,
  at={(11.80cm,0.00cm)},
  anchor=north west,
  width=2.4cm,
  height=2.4cm,
  scale only axis,
  axis equal image,
  hide axis,
  xmin=-2,xmax=2,ymin=-2,ymax=2,
  view={0}{90},
  colormap name=wuemap,
  point meta min=1.055,
  point meta max=2.155,
  title={$b=20$, $L=200$},
  title style={font=\footnotesize,yshift=-0.15cm},
  colorbar,
  colorbar style={width=0.18cm,font=\scriptsize,ytick={1.2,1.6,2.0},height=2.40cm,at={(parent axis.north east)},anchor=north west,xshift=0.25cm}]
\clip (axis cs:0,0) circle [radius=2];
\addplot3[surf,shader=interp,mesh/ordering=x varies,mesh/cols=64] table {figures/tikz/data/rec_b20_L200.dat};
\end{axis}
\draw[black,thick] (recb20L200.center) circle [radius=1.200cm];
\node[below,font=\scriptsize,align=center] at (recb20L200.south) {rel.\ err.\ 0.053\\ time 2\,min};
\draw[gray!50] (2.75cm,0.3cm) -- (2.75cm,-3.25cm);
\end{tikzpicture}
 }
    \caption{Reconstruction of the classical Landweber method and selected reconstructions of the batched Landweber-Kaczmarz method with $b\in\{1,20\}$ and $L\in\{1,200\}$. Below each reconstruction the relative $L^2$-error and the computation time are given.}
    \label{fig:new_num_selected_reconstructions}
\end{figure}

In Figure \ref{fig:new_num_error_results} the total results are plotted and Figure \ref{fig:new_num_selected_reconstructions} shows selected reconstructions.

The computation speed is measured by the number of PDE-solves. For $p$ number of functions and batch size $b$ the number of PDE-solves per iteration is given by
\begin{equation*}
    n_{\textrm{PDE}}= \begin{cases}
        2b \qquad \text{ if the \textit{batched discrepancy} is used}\\
        b+p \quad \text{if the \textit{classical discrepancy }is used}.
    \end{cases}
\end{equation*} 
Note, that the largest part of the computational effort of Algorithm \ref{kaczmarz_algo_eit} lies in solving the PDEs.\\
Figure \ref{fig:new_num_error_results} shows that the quality and computational effort of the reconstructions using the batched discrepancy principle depend significantly on the parameters $L$ and $b$. For a larger batch size $b$, the variance decreases, but at the same time more PDE-solves are required. For a larger aggregation length $L$ the iteration terminates later and hence yields better results.  An exception are the the light blue points with $n_{\text{PDE}}\approx10^5$ and $relative\ error \approx 0.08$. Here one observes that in some runs the stopping criterion fails, terminating only when leaving the range of definition $\sigmad_n\notin \DefSigma$. This is a purely numerical issue; in total the results show: \emph{For a proper choice of parameters $L$ and $b$ the batched discrepancy principle yields equally good approximations as the classical discrepancy principle at significantly lower computation cost}. In our problem the best choice is batch size $b=1$ and aggregation length $L\approx150$.\\
The same effects can be observed in Figure \ref{fig:new_num_selected_reconstructions} showing the selected reconstructions compared with the classical Landweber method.
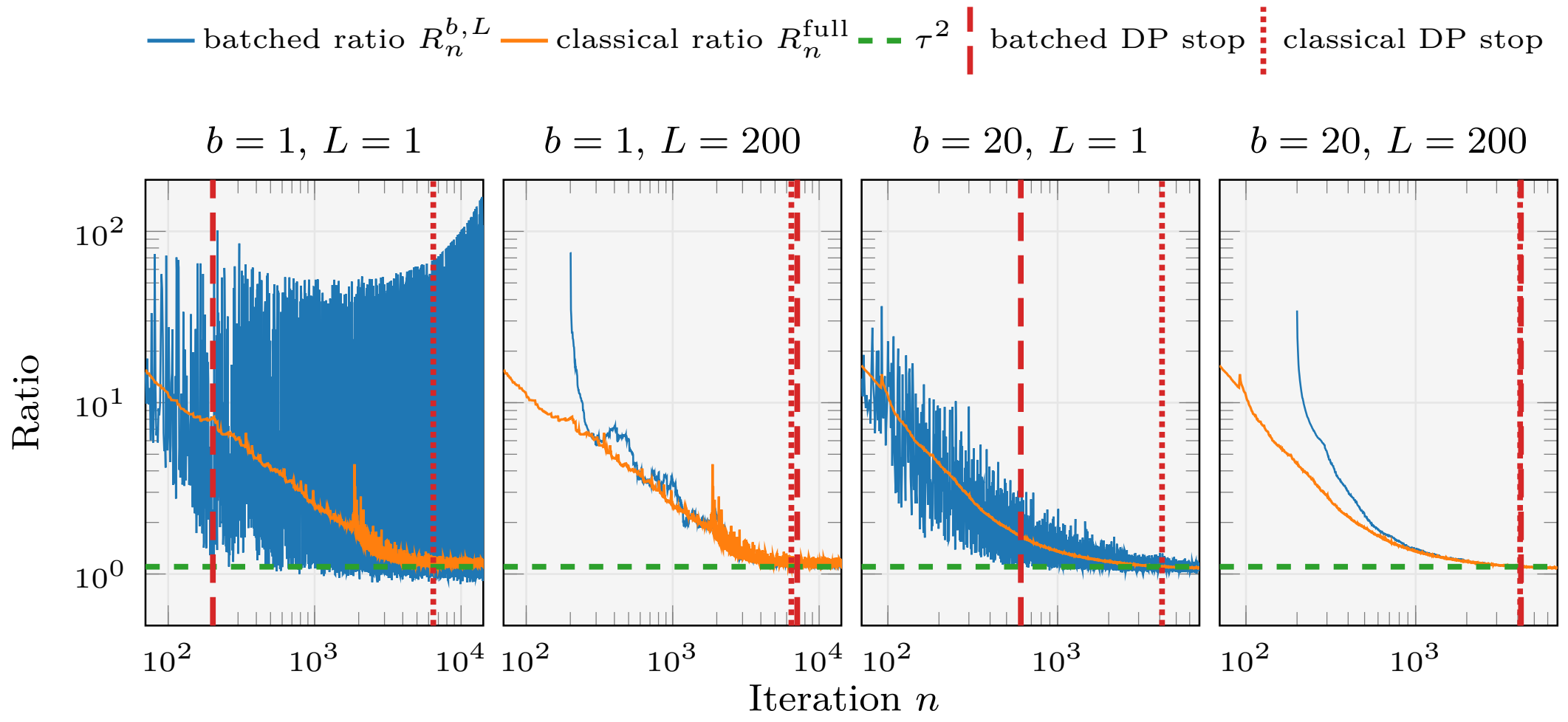
\begin{figure}[ht]
    \centering\tikzset{external/export= true}
    \resizebox{\linewidth}{!}{\definecolor{discBatch}{rgb}{0.1216,0.4667,0.7059}
\definecolor{discFull}{rgb}{1.0000,0.4980,0.0549}
\definecolor{discTau}{rgb}{0.1725,0.6275,0.1725}
\definecolor{discStop}{rgb}{0.8392,0.1529,0.1569}
 \begin{tabular}{@{}c@{}}{\tiny\hbox to 10.35cm{\tikz[baseline=-0.6ex]{\draw[discBatch,line width=0.8pt] (0,0) -- (0.35,0);}\,batched ratio $\DiscRatio^{b,L}_{n}$\hfill
\tikz[baseline=-0.6ex]{\draw[discFull,line width=0.8pt] (0,0) -- (0.35,0);}\,classical ratio $\DiscRatio_n^\text{full}$\hfill
\tikz[baseline=-0.6ex]{\draw[discTau,dashed,very thick] (0,0) -- (0.35,0);}\,$\tau^2$\hfill
\tikz[baseline=-0.6ex]{\path (-0.1,0) -- (0.1,0); \draw[discStop,dash pattern=on 6pt off 3pt,very thick] (0,-0.25) -- (0,0.25);}\,batched DP stop\hfill
\tikz[baseline=-0.6ex]{\path (-0.1,0) -- (0.1,0); \draw[discStop,densely dotted,very thick] (0,-0.25) -- (0,0.25);}\,classical DP stop}}
 \\[0.8em]
\begin{tikzpicture}
\begin{axis}[width=2.5cm,height=3.3cm,scale only axis,xmode=log,ymode=log,xmin=70,ymin=0.5,ymax=200,ytick={1,10,100},xtick={100,1000,10000},tick label style={font=\tiny},label style={font=\scriptsize},title style={font=\scriptsize,yshift=-0.2cm},grid=major,grid style={gray!20},axis background/.style={fill=gray!8},xmax=14183,at={(0.00cm,0cm)},anchor=north west,title={$b=1$, $L=1$},ylabel={Ratio}]
\addplot[discBatch,line width=0.4pt] table[x=k,y=ratio] {figures/tikz/data/disc_b1_L1.dat};
\addplot[discFull,line width=0.5pt] table[x=k,y=ratio] {figures/tikz/data/disc_b1_full.dat};
\addplot[discTau,dashed,very thick,domain=70:14183,samples=2] {1.1025};
\draw[discStop,densely dotted,very thick] (axis cs:6463,0.5) -- (axis cs:6463,200);
\draw[discStop,dash pattern=on 6pt off 3pt,very thick] (axis cs:202,0.5) -- (axis cs:202,200);
\end{axis}
\begin{axis}[width=2.5cm,height=3.3cm,scale only axis,xmode=log,ymode=log,xmin=70,ymin=0.5,ymax=200,ytick={1,10,100},xtick={100,1000,10000},tick label style={font=\tiny},label style={font=\scriptsize},title style={font=\scriptsize,yshift=-0.2cm},grid=major,grid style={gray!20},axis background/.style={fill=gray!8},xmax=14183,at={(2.65cm,0cm)},anchor=north west,title={$b=1$, $L=200$},yticklabels={}]
\addplot[discBatch,line width=0.4pt] table[x=k,y=ratio] {figures/tikz/data/disc_b1_L200.dat};
\addplot[discFull,line width=0.5pt] table[x=k,y=ratio] {figures/tikz/data/disc_b1_full.dat};
\addplot[discTau,dashed,very thick,domain=70:14183,samples=2] {1.1025};
\draw[discStop,densely dotted,very thick] (axis cs:6463,0.5) -- (axis cs:6463,200);
\draw[discStop,dash pattern=on 6pt off 3pt,very thick] (axis cs:7089,0.5) -- (axis cs:7089,200);
\end{axis}
\begin{axis}[width=2.5cm,height=3.3cm,scale only axis,xmode=log,ymode=log,xmin=70,ymin=0.5,ymax=200,ytick={1,10,100},xtick={100,1000,10000},tick label style={font=\tiny},label style={font=\scriptsize},title style={font=\scriptsize,yshift=-0.2cm},grid=major,grid style={gray!20},axis background/.style={fill=gray!8},xmax=6849,at={(5.30cm,0cm)},anchor=north west,title={$b=20$, $L=1$},yticklabels={}]
\addplot[discBatch,line width=0.4pt] table[x=k,y=ratio] {figures/tikz/data/disc_b20_L1.dat};
\addplot[discFull,line width=0.5pt] table[x=k,y=ratio] {figures/tikz/data/disc_b20_full.dat};
\addplot[discTau,dashed,very thick,domain=70:6849,samples=2] {1.1025};
\draw[discStop,densely dotted,very thick] (axis cs:4125,0.5) -- (axis cs:4125,200);
\draw[discStop,dash pattern=on 6pt off 3pt,very thick] (axis cs:607,0.5) -- (axis cs:607,200);
\end{axis}
\begin{axis}[width=2.5cm,height=3.3cm,scale only axis,xmode=log,ymode=log,xmin=70,ymin=0.5,ymax=200,ytick={1,10,100},xtick={100,1000,10000},tick label style={font=\tiny},label style={font=\scriptsize},title style={font=\scriptsize,yshift=-0.2cm},grid=major,grid style={gray!20},axis background/.style={fill=gray!8},xmax=6849,at={(7.95cm,0cm)},anchor=north west,title={$b=20$, $L=200$},yticklabels={}]
\addplot[discBatch,line width=0.4pt] table[x=k,y=ratio] {figures/tikz/data/disc_b20_L200.dat};
\addplot[discFull,line width=0.5pt] table[x=k,y=ratio] {figures/tikz/data/disc_b20_full.dat};
\addplot[discTau,dashed,very thick,domain=70:6849,samples=2] {1.1025};
\draw[discStop,densely dotted,very thick] (axis cs:4125,0.5) -- (axis cs:4125,200);
\draw[discStop,dash pattern=on 6pt off 3pt,very thick] (axis cs:4172,0.5) -- (axis cs:4172,200);
\end{axis}
\node[font=\scriptsize,anchor=north] at (5.17cm,-3.60cm) {Iteration $n$};
\pgfpointanchor{current bounding box}{west}\pgfgetlastxy{\discbbx}{\discbby}
\path ({10.35cm-\discbbx},0);
\end{tikzpicture}
 \end{tabular}
 }
    \caption{Batched and  classical discrepancy ratio $R_{n}^{b,L},R_n^\text{full}$ over the iteration index $n$ for $b\in\{1,20\}$ and $L\in\{1,200\}$. The vertical lines mark the stopping indices of the batched and the classical discrepancy principle.}
    \label{fig:new_num_discrepancy_ratios}
\end{figure}
To explain the reason behind the influence of the parameters on the results we plotted the batched and classical discrepancy ratio over the iteration index $n$ for $b\in\{1,20\}$ and $L\in\{1,200\}$ (Figure \ref{fig:new_num_discrepancy_ratios}). Here we observe, that for batch size $b=20$ the ratio is less scattered than for batch size $b=1$. In the plots with $L=1$, one can also see that, as explained in Section \ref{section:batch_DP}, the batched discrepancy for $L=1$ is an unbiased estimator of the classical discrepancy. However, in this case the batched ratio has a large variance, hence the approximation of the classical ratio is unstable and the stopping too early. The problem is solved choosing larger values of $L$, which gives a reduced variance, but an additional error occurs due to the consideration of previous iterates. The plots of $L=200$ show, that the additional error is limited to the first iterates. For later iterates the batched ratio yields a good approximation of the classical ratio and therefore similar stopping at lower computational cost.

\section{Concluding Remarks}
In this work, we analyzed a \emph{stochastic mini-batch Landweber-Kaczmarz method} with a new a posteriori stopping: the batched discrepancy principle. Under the tangential cone condition we proved the regularization properties for a priori stopping and the introduced batched discrepancy principle. Numerical experiments for the electrical impedance tomography problem show that: For a proper choice of the aggregation parameter $L$ and the batch size $b$, the batched discrepancy principle yields equally good approximations as the classical discrepancy principle at significantly lower computational cost. Building on our work, there are some questions that might be considered for future research:
\begin{enumerate}
    \item Can the method be generalized to a Banach space setup, i.e. the parameter space $X$ and/or data spaces $Y_i$ are Banach spaces? This will cover inverse problems with nonsmooth solutions.
    \item What convergence rates can be proven?
    \item How does the method extend to accelerations, like for example the Nesterov acceleration?
\end{enumerate}
\section{Acknowledgment}
Rommel R. Real is supported by Centre International de Math\'{e}matiques Pures et Appliqu\'{e}es (CIMPA) and the Abdus Salam International Centre for Theoretical Physics (ICTP) through the CIMPA-ICTP Research in Pairs fellowship programme.

\section{Disclaimer}
During the writing of this paper, artificial intelligence was used for graphic design and proof-reading parts of the manuscript. We carefully checked all AI suggested content and take full responsibility for all mathematical content, software, numerical results, interpretations and conclusions.

\clearpage
\bibliographystyle{abbrv}

\begin{thebibliography}{10}

\bibitem{Bahouri_2011}
H.~Bahouri, J.-Y. Chemin, and R.~Danchin.
\newblock {\em Fourier analysis and nonlinear partial differential equations},
  volume 343 of {\em Grundlehren der mathematischen Wissenschaften [Fundamental
  Principles of Mathematical Sciences]}.
\newblock Springer, Heidelberg, 2011.

\bibitem{Bhattacharya2016}
R.~Bhattacharya and E.~C. Waymire.
\newblock {\em A basic course in probability theory}.
\newblock Universitext. Springer, Cham, second edition, 2016.

\bibitem{BlanchardHoffmannReiss2018a}
G.~Blanchard, M.~Hoffmann, and M.~Rei\ss.
\newblock Early stopping for statistical inverse problems via truncated {SVD}
  estimation.
\newblock {\em Electron. J. Stat.}, 12(2):3204--3231, 2018.

\bibitem{BlanchardHoffmannReiss2018b}
G.~Blanchard, M.~Hoffmann, and M.~Rei\ss.
\newblock Optimal adaptation for early stopping in statistical inverse
  problems.
\newblock {\em SIAM/ASA J. Uncertain. Quantif.}, 6(3):1043--1075, 2018.

\bibitem{Borcea_2002}
L.~Borcea.
\newblock Electrical impedance tomography.
\newblock {\em Inverse Problems}, 18(6):R99--R136, 2002.

\bibitem{Bottou2010}
L.~Bottou.
\newblock Large-scale machine learning with stochastic gradient descent.
\newblock In {\em Proceedings of {COMPSTAT}'2010}, pages 177--186.
  Physica-Verlag/Springer, Heidelberg, 2010.

\bibitem{Bottou2018}
L.~Bottou, F.~E. Curtis, and J.~Nocedal.
\newblock Optimization methods for large-scale machine learning.
\newblock {\em SIAM Rev.}, 60(2):223--311, 2018.

\bibitem{Chen2018}
K.~Chen, Q.~Li, and J.-G. Liu.
\newblock Online learning in optical tomography: a stochastic approach.
\newblock {\em Inverse Problems}, 34(7):075010, 2018.

\bibitem{Clason2019}
C.~Clason and V.~H. Nhu.
\newblock Bouligand-{L}andweber iteration for a non-smooth ill-posed problem.
\newblock {\em Numer. Math.}, 142(4):789--832, 2019.

\bibitem{de2008steepest}
A.~De~Cezaro, M.~Haltmeier, A.~Leit\~{a}o, and O.~Scherzer.
\newblock On steepest-descent-{K}aczmarz methods for regularizing systems of
  nonlinear ill-posed equations.
\newblock {\em Appl. Math. Comput.}, 202(2):596--607, 2008.

\bibitem{Defazio2014}
A.~Defazio, F.~Bach, and S.~Lacoste-Julien.
\newblock {SAGA}: A fast incremental gradient method with support for
  non-strongly convex composite objectives.
\newblock In {\em Advances in Neural Information Processing Systems},
  volume~27. Curran Associates, Inc., 2014.

\bibitem{Engl2000}
H.~W. Engl, M.~Hanke, and A.~Neubauer.
\newblock {\em Regularization of inverse problems}, volume 375 of {\em
  Mathematics and its Applications}.
\newblock Kluwer Academic Publishers Group, Dordrecht, 1996.

\bibitem{evans2010partial}
L.~C. Evans.
\newblock {\em Partial differential equations}, volume~19 of {\em Graduate
  Studies in Mathematics}.
\newblock American Mathematical Society, Providence, RI, second edition, 2010.

\bibitem{Gitman2019}
I.~Gitman, H.~Lang, P.~Zhang, and L.~Xiao.
\newblock Understanding the role of momentum in stochastic gradient methods.
\newblock In {\em Advances in Neural Information Processing Systems},
  volume~32. Curran Associates, Inc., 2019.

\bibitem{Gu2025}
R.~Gu, Z.~Fu, B.~Han, and H.~Fu.
\newblock Stochastic gradient descent method with convex penalty for ill-posed
  problems in {B}anach spaces.
\newblock {\em Inverse Problems}, 41(5):055003, 2025.

\bibitem{Haltmeier2007}
M.~Haltmeier, A.~Leit\~{a}o, and O.~Scherzer.
\newblock Kaczmarz methods for regularizing nonlinear ill-posed equations. {I}.
  {C}onvergence analysis.
\newblock {\em Inverse Probl. Imaging}, 1(2):289--298, 2007.

\bibitem{Hanke_1991}
M.~Hanke.
\newblock Accelerated {L}andweber iterations for the solution of ill-posed
  equations.
\newblock {\em Numer. Math.}, 60(3):341--373, 1991.

\bibitem{Hanke1995}
M.~Hanke, A.~Neubauer, and O.~Scherzer.
\newblock A convergence analysis of the {L}andweber iteration for nonlinear
  ill-posed problems.
\newblock {\em Numer. Math.}, 72(1):21--37, 1995.

\bibitem{Herman1993}
G.~T. Herman and L.~B. Meyer.
\newblock Algebraic reconstruction techniques can be made computationally
  efficient.
\newblock {\em IEEE Trans. Med. Imaging}, 12(3):600--609, 1993.

\bibitem{Hertrich2022}
J.~Hertrich and G.~Steidl.
\newblock Inertial stochastic {PALM} and applications in machine learning.
\newblock {\em Sampl. Theory Signal Process. Data Anal.}, 20(1):4, 2022.

\bibitem{Huang2025}
J.~Huang, Q.~Jin, X.~Lu, and L.~Zhang.
\newblock On early stopping of stochastic mirror descent method for ill-posed
  inverse problems.
\newblock {\em Numer. Math.}, 157(2):539--571, 2025.

\bibitem{HuckerReiss2025}
L.~Hucker and M.~Rei\ss.
\newblock Early stopping for conjugate gradients in statistical inverse
  problems.
\newblock {\em Numer. Math.}, 157(5):1739--1791, 2025.

\bibitem{Jahn2020}
T.~Jahn and B.~Jin.
\newblock On the discrepancy principle for stochastic gradient descent.
\newblock {\em Inverse Problems}, 36(9):095009, 2020.

\bibitem{Jia2010}
N.~Jia and E.~Y. Lam.
\newblock Machine learning for inverse lithography: using stochastic gradient
  descent for robust photomask synthesis.
\newblock {\em J. Opt.}, 12(4):045601, 2010.

\bibitem{JinKereta2023}
B.~Jin and {\v{Z}}.~Kereta.
\newblock On the convergence of stochastic gradient descent for linear inverse
  problems in {B}anach spaces.
\newblock {\em SIAM J. Imaging Sci.}, 16(2):671--705, 2023.

\bibitem{Jin2026a}
B.~Jin, {\v{Z}}.~Kereta, and Y.~Xia.
\newblock Stochastic gradient descent for nonlinear inverse problems in
  {B}anach spaces.
\newblock {\em SIAM J. Imaging Sci.}, 19(2):1244--1271, 2026.

\bibitem{JinLu2019}
B.~Jin and X.~Lu.
\newblock On the regularizing property of stochastic gradient descent.
\newblock {\em Inverse Problems}, 35(1):015004, 2019.

\bibitem{Jin2025c}
B.~Jin and Z.~Zhou.
\newblock On the convergence of stochastic variance reduced gradient for linear
  inverse problems.
\newblock {\em Inverse Problems}, 42(4):045006, 2026.

\bibitem{Jin2020}
B.~Jin, Z.~Zhou, and J.~Zou.
\newblock On the convergence of stochastic gradient descent for nonlinear
  ill-posed problems.
\newblock {\em SIAM J. Optim.}, 30(2):1421--1450, 2020.

\bibitem{Jin_2022SVRG}
B.~Jin, Z.~Zhou, and J.~Zou.
\newblock An analysis of stochastic variance reduced gradient for linear
  inverse problems.
\newblock {\em Inverse Problems}, 38(2):025009, 2022.

\bibitem{Jin2025b}
Q.~Jin and L.~Chen.
\newblock Stochastic variance reduced gradient method for linear ill-posed
  inverse problems.
\newblock {\em Inverse Problems}, 41(5):055014, 2025.

\bibitem{Jin2023b}
Q.~Jin, X.~Lu, and L.~Zhang.
\newblock Stochastic mirror descent method for linear ill-posed problems in
  {B}anach spaces.
\newblock {\em Inverse Problems}, 39(6):065010, 2023.

\bibitem{Jin2013}
Q.~Jin and W.~Wang.
\newblock Landweber iteration of {K}aczmarz type with general non-smooth convex
  penalty functionals.
\newblock {\em Inverse Problems}, 29(8):085011, 2013.

\bibitem{Johnson2013}
R.~Johnson and T.~Zhang.
\newblock Accelerating stochastic gradient descent using predictive variance
  reduction.
\newblock In {\em Advances in Neural Information Processing Systems},
  volume~26, pages 315--323. Curran Associates, Inc., 2013.

\bibitem{Kaltenbacher2008}
B.~Kaltenbacher, A.~Neubauer, and O.~Scherzer.
\newblock {\em Iterative regularization methods for nonlinear ill-posed
  problems}, volume~6 of {\em Radon Series on Computational and Applied
  Mathematics}.
\newblock Walter de Gruyter, Berlin, 2008.

\bibitem{Kindermann_2022}
S.~Kindermann.
\newblock On the tangential cone condition for electrical impedance tomography.
\newblock {\em Electron. Trans. Numer. Anal.}, 57:17--34, 2022.

\bibitem{Kowar_2002}
R.~Kowar and O.~Scherzer.
\newblock Convergence analysis of a {L}andweber-{K}aczmarz method for solving
  nonlinear ill-posed problems.
\newblock In V.~G. Romanov, S.~I. Kabanikhin, Y.~E. Anikonov, and A.~L.
  Bukhgeim, editors, {\em Ill-posed and inverse problems}, pages 253--270. VSP,
  Zeist, 2002.

\bibitem{Leitao2016}
A.~Leit\~{a}o and B.~F. Svaiter.
\newblock On projective {L}andweber-{K}aczmarz methods for solving systems of
  nonlinear ill-posed equations.
\newblock {\em Inverse Problems}, 32(2):025004, 2016.

\bibitem{Li2018}
H.~Li and M.~Haltmeier.
\newblock The averaged {K}aczmarz iteration for solving inverse problems.
\newblock {\em SIAM J. Imaging Sci.}, 11(1):618--642, 2018.

\bibitem{needell2014stochastic}
D.~Needell, N.~Srebro, and R.~Ward.
\newblock Stochastic gradient descent, weighted sampling, and the randomized
  {K}aczmarz algorithm.
\newblock {\em Math. Program.}, 155(1--2):549--573, 2016.

\bibitem{Nguyen2017}
L.~M. Nguyen, J.~Liu, K.~Scheinberg, and M.~Tak{\'a}{\v{c}}.
\newblock {SARAH}: A novel method for machine learning problems using
  stochastic recursive gradient.
\newblock In {\em Proceedings of the 34th International Conference on Machine
  Learning}, volume~70 of {\em Proc. Mach. Learn. Res.}, pages 2613--2621.
  PMLR, 2017.

\bibitem{Robbins1951}
H.~Robbins and S.~Monro.
\newblock A stochastic approximation method.
\newblock {\em Ann. Math. Statistics}, 22:400--407, 1951.

\bibitem{schmidt2017minimizing}
M.~Schmidt, N.~Le~Roux, and F.~Bach.
\newblock Minimizing finite sums with the stochastic average gradient.
\newblock {\em Math. Program.}, 162(1--2):83--112, 2017.

\bibitem{ShalevShwartz2013}
S.~Shalev-Shwartz and T.~Zhang.
\newblock Stochastic dual coordinate ascent methods for regularized loss
  minimization.
\newblock {\em J. Mach. Learn. Res.}, 14:567--599, 2013.

\bibitem{Strohmer2009}
T.~Strohmer and R.~Vershynin.
\newblock A randomized {K}aczmarz algorithm with exponential convergence.
\newblock {\em J. Fourier Anal. Appl.}, 15(2):262--278, 2009.

\bibitem{Zhang2013}
L.~Zhang, M.~Mahdavi, and R.~Jin.
\newblock Linear convergence with condition number independent access of full
  gradients.
\newblock In {\em Advances in Neural Information Processing Systems},
  volume~26. Curran Associates, Inc., 2013.

\bibitem{Zhang2015}
Y.~Zhang and L.~Xiao.
\newblock Stochastic primal-dual coordinate method for regularized empirical
  risk minimization.
\newblock In {\em Proceedings of the 32nd International Conference on Machine
  Learning}, volume~37 of {\em Proc. Mach. Learn. Res.}, pages 353--361. PMLR,
  2015.

\end{thebibliography}

\end{document}